\documentclass[preprint,3p,11pt,authoryear]{elsarticle}

\usepackage{amssymb}
\usepackage{amsmath}
\usepackage{multirow}

\journal{European Journal of Operational Research}

\newcommand{\bA}{ \mathbf{A} }
\newcommand{\ba}{ \mathbf{a} }
\newcommand{\bB}{ \mathbf{B} }
\newcommand{\bb}{ \mathbf{b} }

\newcommand{\bc}{ \mathbf{c} }
\newcommand{\bD}{ \mathbf{D} }

\newcommand{\bF}{ \mathbf{F} }
\newcommand{\bE}{ \mathbf{E} }

\newcommand{\bu}{ \mathbf{u} }

\newcommand{\bW}{ \mathbf{W} }

\newcommand{\bX}{ \mathbf{X} }
\newcommand{\bx}{ \mathbf{x} }

\newcommand{\by}{ \mathbf{y} }

\newcommand{\FO}{ \mathbf{FO} }
\newcommand{\IL}{ \text{IL} }
\newcommand{\ILL}{ \text{IL} }
\newcommand{\bepsilon}{ \mathbf{\epsilon} }

\newcommand{\bzero}{\mathbf{0}}

\newcommand{\cD}{\mathcal{D}}
\newcommand{\cJ}{\mathcal{J}}

\newcommand{\IPP}{\text{RTP}}

\newcommand{\cS}{\mathcal{S}}
\newcommand{\MIL}{\text{MIL}}
\newcommand{\todo}[1]{\textit{\textcolor{red}{#1}}}

\usepackage[nameinlink]{cleveref}
\crefname{assumption}{Assumption}{Assumptions}
\crefname{lemma}{Lemma}{Lemmas}
\crefname{theorem}{Theorem}{Theorems}
\crefname{corollary}{Corollary}{Corollaries}
\crefname{proposition}{Proposition}{Propositions}
\crefname{claim}{Claim}{Claims}
\crefname{definition}{Definition}{Definitions}
\crefname{subclaim}{Subclaim}{Subclaims}
\crefname{procedure}{Procedure}{Procedures}
\crefname{algorithm}{Algorithm}{Algorithms}
\crefname{example}{Example}{Examples}
\crefname{figure}{Figure}{Figures}
\crefname{section}{Section}{Sections}
\crefname{appendix}{Appendix}{Appendices}
\crefname{table}{Table}{Tables}
\crefname{equation}{}{}

\usepackage{makecell}

\newtheorem{thm}{Theorem}

\newtheorem{prop}{Proposition}
\newdefinition{rmk}{Remark}
\newproof{pf}{Proof}
\newproof{pot}{Proof of Theorem \ref{thm2}}

\usepackage{array}

\begin{document}

\begin{frontmatter}



\title{Improving Observed Decisions' Quality using Inverse Optimization: A Radiation Therapy Treatment Planning Application
} 

\author[1]{Farzin Ahmadi}
{\ead{fahmadi1@jhu.edu}}
\author[2]{Todd R. McNutt}
\ead{tmcnutt1@jhmi.edu}
\author[1]{Kimia Ghobadi\corref{cor1}}
\ead{kimia@jhu.edu}
\cortext[cor1]{Corresponding author}
\affiliation[1]{organization={Department of Civil and Systems Engineering and  Malone Center for Engineering in Healthcare, Johns Hopkins University},
city={Baltimore},
country={USA}}
\affiliation[2]{organization={Department of Radiation Oncology and Molecular Radiation Sciences, Johns Hopkins School of Medicine},
city={Baltimore},
country={USA}}

\begin{abstract}
In many applied optimization settings, parameters that define the constraints may not guarantee the best possible solution, and superior solutions might exist that are infeasible for the given parameter values.  Removing such constraints, re-optimizing, and evaluating the new solution may be insufficient, as the optimizer's preferences in selecting the existing solutions might be lost. To address this issue, we present an inverse optimization-based model that takes an observed solution as input and aims to improve upon it by projecting onto desired hyperplanes or expanding the feasible set while balancing the distance to the observed decision to preserve the optimizer's preferences.

We demonstrate the applicability of the model in the context of radiation therapy treatment planning, an essential component of cancer treatment. Radiation therapy treatment planning is typically guided by expert-driven guidelines that define the optimization problem but remain mostly general. Our model provides an automated framework that learns new plans from available plans based on given clinical criteria, optimizing the desired effect without compromising the remaining constraints.

The proposed approach is applied to a cohort of four prostate cancer patients, and the results demonstrate improvements in dose-volume histograms while maintaining comparable target coverage to clinically acceptable plans. By optimizing the parameters of the treatment planning problem and exploring the Pareto frontier, our methodology uncovers previously unattainable solutions that enhance organ-at-risk sparing without sacrificing target coverage. The framework's ability to handle multiple organs-at-risk and various dose-volume constraints highlights its flexibility and potential for application to diverse radiation therapy treatment planning scenarios.

\end{abstract}

\end{frontmatter}



\section{Introduction} \label{sec:Introduction}


Radiation therapy is a widely utilized cancer treatment modality, with approximately 50\% of cancer patients receiving radiation therapy and 30\% of cancer survivors having undergone this form of treatment (Miller et al., 2016). However, optimizing radiation therapy treatment plans remains a significant challenge due to the need to balance high radiation doses to target structures while minimizing exposure to healthy tissues. This is crucial for preventing potential damage to normal tissues and reducing long-term side effects. Moreover, clinical objectives for radiation therapy plans are often broad and not personalized, and plan improvements are usually possible. To address these challenges, we aim to develop an automated treatment improvement strategy that modifies optimization parameters to provide better plans with mathematical guarantees using clinically acceptable plans as input. Our focus is primarily on Intensity Modulated Radiation Therapy (IMRT), a method that allows for precise conformity of radiation to target structures \citep{taylor2004intensity}, and its application in prostate cancer treatment planning.

\begin{figure}[htbp]
\begin{center}
\includegraphics[width =0.8 \linewidth]{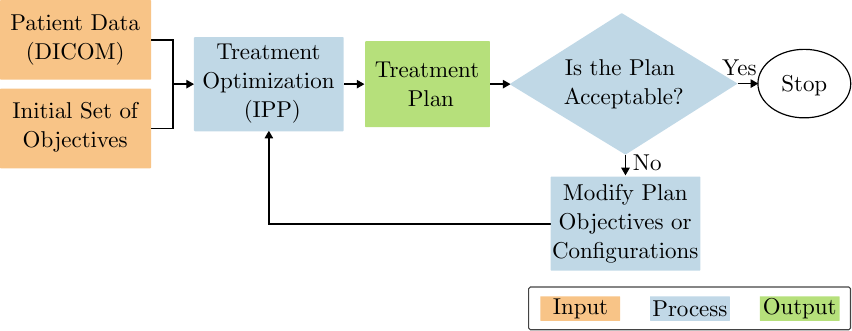}
\caption{\footnotesize Existing treatment planning frameworks input patient data (usually stored in DICOM standard files) and an initial set of objectives. The process requires potentially numerous rounds of manual modifications to plan parameters or configurations until a clinically acceptable plan is achieved. This iterative approach can be time-consuming and may not guarantee optimality in the resulting treatment plans.} \label{Fig:IPP_framework}
\end{center}
\end{figure}

Current radiation therapy treatment planning often requires many rounds of optimization \citep{wu2009adaptive}, as dose-volume  (DV) objectives drive the optimization and parameters are progressively set by a treatment planner until a clinically acceptable plan is achieved \citep{wang2013quality} (as shown in \ref{Fig:IPP_framework}). This process remains heavily expert-driven and time-consuming, with plan qualities varying among different experts and planners \citep{wu2014improved,petit2012increased,chung2008can}. Previous studies have attempted to develop data-driven models to provide better objectives for future patients based on quantitative comparisons between features of new and previous patients \citep{mahmood2018automated, sahiner2019deep, wang2013quality, petit2012increased}. For instance, machine learning applications in radiation therapy have been detailed in \cite{el2015machine} and  \cite{kang2015machine}. While efforts have been made to provide user-friendly environments for experts to compare different plans, Studies aiming to automatically improve the available plans in terms of sparing organs at risk (OARs) and target coverage through mathematical modeling are rare.

Optimization models play a crucial role in characterizing the trade-offs between different objectives and constraints, and in identifying the limitations imposed by the current problem formulation. In many optimization problems, constraints are introduced to model the requirements and limitations of the system under consideration. These requirement constraints, which define the feasible region of the optimization problem, are typically determined by the problem parameters and/or expert opinions. However, in practice, the parameter values used to specify these constraints may not always result in the best possible solution. In fact, there may exist better solutions that are rendered infeasible by the current parameter settings \citep{gorissen2013mixed}. Simply removing such constraints, re-optimizing, and evaluating the new solution might be insufficient, as the optimizer's preferences in choosing the existing solutions might disappear. To address this issue, in this work, we present an inverse optimization-based model that inputs an observed solution and aims to improve upon it by projecting onto desired hyperplanes or expanding the feasible set while balancing the distance to the observed decision to preserve optimizer preferences. Inverse optimization has been increasingly applied in various fields, including healthcare \citep{erkin2010eliciting, ghobadi2018robust, ahmadi2020inverse}, finance \citep{bertsimas2012inverse}, and transportation \citep{chow2012inverse}. In radiation therapy, inverse optimization has been used to learn unknown parameters of the optimization problem \citep{babier2018inverse,chan2014generalized,chan2020inverse}. Other streams of literature have attempted at posing the treatment planning optimization problem as a tractable optimization problem include posing linear Programming models \citep{romeijn2006new,holder2005radiotherapy,morrill1991dose,craft2007local} or quadratic optimization models. More recent  works have considered methods for learning objectives and objective weights for radiation therapy treatment plans \citep{shen2020operating,boutilier2015models,ge2019knowledge}. However, investigations into the mathematical structure of the radiation therapy treatment planning optimization problem resulting from preset expert-driven objectives have been limited.

In this work, we employ inverse optimization-based methods to improve existing plans of a patient, which we assume carry some notions of the decision maker preferences. We propose improved plans in terms of DV objectives without compromising target or other OAR goals, if such plans exist.  By optimizing some of the parameters of the treatment planning optimization problem while providing a solution that adheres to the improvements, our methodology differs from previous data-driven models that are considered for recovering or generating personalized objectives using existing data for future patients. Naturally, employing all potential directions and providing improved plans will allow experts to choose from a set of optimal plans for future patients. We additionally show how employing such a methodology will allow alternative characterizations of the Pareto frontier of the optimization problem and how the proposed new objectives/constraints shape the Pareto frontier. The main contributions of this paper are as follows:

1. We develop an inverse optimization-based model that learns new and improved radiation therapy treatment plans from observed plans while preserving optimizer preferences.

2. We demonstrate the applicability of the model in the setting of prostate cancer radiation therapy treatment planning and show improvements in dose-volume histograms with comparable target coverage to clinically acceptable plans.

3. We explore the implications of our methodology on the Pareto frontier of the multi-objective representation of the treatment planning problem before and after improving a single dose limit.

The remainder of this paper is organized as follows. Section \ref{sec:IL_methodology} details the key concepts and methodologies, while Section \ref{sec:RT_methodology} applies the developed methodologies to the plan improvement application. In Section \ref{sec:results}, we validate the models by providing the results of employing the methodology on four prostate cancer patients. Discussions on the applicability and practicality of the methods introduced in this work and conclusions are presented in Section \ref{sec:discussions}.

\section{Discovering Optimal Solutions From Current Decisions}\label{sec:IL_methodology}
In this section, we characterize the setting of automating the discovery of optimal solutions from observed current decisions. We pose a general description of the problem and attempt to model the problem for the specific case where the underlying optimization problem is linear. We explore the properties of the model for the linear case and connect the model to the inverse optimization literature. 

\subsection{General Form of the Problem} \label{sec:IL}
Let $\FO (\bu, \theta) = \left\{ \max_{\bx} f(\bx, \bu, \theta)| \bx \in \Omega(\bu, \theta) \right\}$ be an optimization problem in general form where $\bx$ is the decision variable and $f$ is the objective as a function of the variables $\bx$ and the parameters $(\bu, \theta)$, and $\Omega(\bu, \theta)$ is the feasible set for $\FO$. Since optimal values of $\bx$ can be found by solving $\FO$ when all parameters are available, we assume that $\bu \in \mathbb{R}^{m_1}$ comprises of missing or unknown parameters and $\theta \in \mathbb{R}^{m_2}$ represent known parameters. Hence, we do not readily have access to optimal values of $\bx$. 

In the case where $\{\bu\} \neq \emptyset$, inverse optimization \citep{Ahuja2001, aswani2018inverse} aims to recover unknown parameters of $\FO$ when solution(s) of the problem are observed. When focused on learning the optimal solution of $\FO$, \cite{ahmadi2020inverse} provide an alternative modeling of the inverse optimization problem (denoted as Inverse Learning ($\IL$)) for the linear optimization case. Building on that model, for the case decision(s) $\bX_0 \subset \mathbb{R}^{n}$ are given and an optimal solution $\bx$ is desired, a similar model can be considered for our $\FO$ as follows:
\begin{subequations} \label{IL_general}
\begin{align}
\IL(\theta,\bX_0): 
\underset{\bx, \bu}{\text{minimize}} & \quad \cD(\bx, \bX_0) \label{IL_obj}\\
\text{subject to} 
& \quad \bx \in \Omega^{opt}(\bu,\theta),  \label{IL_Const} \\
& \quad  \bx \in \mathbb{R}^n, \quad  \bu \in \mathbb{R}^m. \label{IL_range}
\end{align}
\end{subequations}

In the above model, $\bx$ and $\bu$ constitute the decision variables. $\Omega^{opt}(\bu,\theta)$ characterizes the set of all solutions that can be optimal for $\FO(\bu, \theta)$. For any feasible solution of $\IL$, a solution $\bx$ to $\FO(\bu, \theta)$ is found such that it is optimal for that realization of the unknown parameter $\bu$, as shown by constraint \eqref{IL_Const}. The optimal solution of $\IL$ minimizes a loss function $\cD$ with regard to $\bX_0$, available and/or observed solution(s) of $\FO$. While in general, characterizing the set of optimal solutions $\Omega^{opt}(\bu, \theta)$ is not trivial, we note that for the most general case for $f$, $\Omega^{opt}(\bu, \theta)$ could be even equal to $\Omega$. However, for the cases where $f$ is monotonic (at least over $\Omega$), generally, $\Omega^{opt}(\bu, \theta) \subseteq \Omega$. Finally, we note that in general, if $f$ is not monotonic, the solution to $\IL$ can be the average point of all observations in $\bX_0$ as it is generally possible to find the combination of parameters that render the average optimal. It is mainly in the case where $f$ is monotonic in structure that the average point of all observations in $\bX_0$ is not guaranteed to be the optimal solution of IL. As such, we focus our attention on cases where $f$ is monotonic, as this case also aligns well with our application. Specifically, we focus on the case where $\FO$ is linear and provide tractable models in the following.

 \subsection{Explicit Formulation for Linear $\FO$}
In our radiation therapy application, we aim to pose the (forward) treatment planning optimization problem as a linear optimization model for simplicity and tractability, and while $\IL$ is a bi-level optimization problem, it can be explicitly formulated by outlining optimality conditions for $\bx$ for a linear $\FO$ \citep{ahmadi2020inverse}. 
In this case, the forward optimization problem can be characterized as $\FO (\bc, \bA, \bb) = \left\{ \max_{\bx} \bc^T \bx| \bA \bx \leq \bb, \bx\in \mathbb{R}^n \right\}$ and $\Omega = \left\{ \bx \in \mathbb{R}^n| \bA \bx \leq \bb \right\}$ is the feasible set and we assume $\Omega$ to be non-empty, full dimensional and free of redundant constraints. For the case where $\Omega$ is known and the cost vector $\bc$ is unknown, $\bu = \bc$ and $\theta = [\bA \ \bb]$. Then, the linear version of $\IL$ is formulated based on duality conditions that ensure $\bx$ is optimal for a realization of $\bc$. The inverse learning model for a linear $\FO$ can then be formulated as follows for the case where $\bu = \bc$.
\begin{subequations} \label{IL_Linear}
\begin{align}
\ILL({\bX_0, \bA,\bb}): \underset{\bc, \bx, \by}{\text{minimize}} & 
\quad \cD(\bx, \bX_0) \\
\text{subject to} 
& \quad \bA \bx \leq \bb,   \label{IOPrimal Feasblility1}\\ 
& \quad \bA' \by  =\bc, \label{IODualFeas1}\\ 
& \quad \bc' \bx = \bb' \by ,   \label{IOStrongDual}\\ 
& \quad \sum_{j \in \cJ} y_j =1,   \label{IORegularization}\\
& \quad \by \ge \bzero,  \label{IODualFeas2}\\
& \quad  \bx \in \mathbb{R}^n, \quad \bc \in \mathbb{R}^n, \quad\by \in \mathbb{R}^m_{\geq 0}. \label{IOrange}
\end{align}
\end{subequations}

In the above model, constraint \eqref{IOPrimal Feasblility1} ensures that the learned solution $\bx$ is feasible for $\FO$, constraints \eqref{IODualFeas1}, \eqref{IODualFeas2} and \eqref{IOStrongDual} provide dual feasibility and strong duality conditions for $(\bx,\by)$. The objective function of the inverse problem is similar to formulation \ref{IL_general}. Note that each feasible solution of $\IL$ consists of a cost vector $\bc$ and a pair of optimal primal-dual solutions for $\FO$. For each such solution, $\bx$ is indeed an optimal solution to $\FO(\bc, \bA, \bb)$. However, $\IL$ finds the solution $(\bc, \bx, \by)$ that minimizes a loss with relation to observed decision(s) $\bX_0$. One such loss function can simply be the distance between the optimal solution to be learned ($\bx$) and the observed decision(s), signalling adherence to the policies that derived the observed decision(s) while also ensuring mathematical optimality of $\bx$ for the realized $\FO$. 

The inverse learning model shown in formulation \eqref{IL_Linear} is based on the assumption that no knowledge is available on the cost vector $\bc$. This assumption, however, is too restrictive in many decision making environments where expert concensus exists with regards to decision goals. In such settings, a notion of the direction of $\bc$ or the sign of each element of $\bc$ is readily available. For instance, in radiation therapy treatment, when outlining the multi-objective optimization problem, while the exact values of different weights are not readily known, one can easily assign signs to different objective terms. In a minimization setting, the objective measuring the mean dose levels of an organ at risk has a positive sign and an objective pertaining to the mean dose of a target structure has a negative sign. Additionally, for dose volume constraints that are prevalent in radiation therapy treatment planning, more stringent dose limits on organs at risk are desirable as they provide better radiation sparing of the organ, provided that the optimization remains feasible. In what follows, we develop a modified $\IL$ model that is capable of improving observed decisions when implicit information on the cost vector is available.

\subsection{Improving Observed Decisions via Inverse Learning}

Let $\bx_0$ be an observed decision for the linear representation of $\FO$. In addition to recovering a cost vector, we assume that a non-empty subset of constraints are desired to be improved. In other words, left-hand-side values are subject to increase as much as possible so that a perturbed solution is found that (a) is optimal for the recovered cost vector and (b) is feasible for all remaining constraints and binding the new realizations of the constraints subject to improvement. For the constraints subject to improvement, we assume the dual variables $y$ are non-zero, ensuring that these constraints are binding for the learned optimal solution. We assume implicit knowledge on the cost vector in the form of dual variables, allowing the models to pursue superior realization of right-hand-side values of select constraints. In this case, the inverse process generally pursues two objectives. \emph{First}, similar to Section \ref{sec:IL}, the inverse process aims at finding an optimal solution that is minimally perturbed from the observed decision $\bx_0$. \emph{Second}, due to the available knowledge of the dual variables $\by$, it is desired to maximize the $\FO$ objective function resulting from the learned optimal solution. In general, these two objectives might be conflicting and warrant a new objective function for the inverse learning problem. Hence, a  general inverse learning formulation can then be characterized as follows. In the following, let $\cJ$ be the set of indices of all constraints, $\hat{\cJ}$ and $\bar{\cJ}$ be the set of indices of desired constraints for improvements, and the rest of the constraints respectively ($\hat{\cJ} \cup \bar{\cJ} = \cJ$) and $\by = \left [\hat{\by}; \,\bar{\by} \right]$. We assume 
$\hat{\cJ} \neq \emptyset$ unless stated otherwise. Finally, let $\bb = \left [\hat{\bb}; \,\bar{\bb} \right]$ be the known and unknown right-hand-side parameters with vectors $\bb_L \in \mathbb{R}^m$ and $\bb_U \in \mathbb{R}^m$ as upper and lower bounds of $\bb$. Additionally, without loss of generality, we always assume access to initial values $\hat{\bb}_0$ that form $\Omega$. As such, we have the following model:
\begin{subequations} \label{IL_General}
\begin{align}
\IL_g({\bx_0, \bA, \bar{\bb},\bar{\by} = 0, \omega}): \underset{\bc, \bx, \hat{\by}, \hat{\bb}}{\text{minimize}} & 
\quad \omega\cD(\bx, \bx_0) - (1 - \omega) (\left [\hat{\bb}; \,\bar{\bb} \right]' \left [\hat{\by}; \,\bar{\by} \right])\\
\text{subject to} 
& \quad \Bar{\bA} \bx \leq \Bar{\bb},   \label{MILPrimal Feasblility1}\\ 
& \quad \hat{\bA} \bx \leq \hat{\bb},   \label{MILPrimal Feasblility2}\\ 
& \quad \bA' \left [\hat{\by}; \,\bar{\by} \right]  =\bc, \label{MILDualFeas1}\\ 
& \quad \bc' \bx = \left [\hat{\bb}; \,\bar{\bb} \right]' \left [\hat{\by}; \,\bar{\by} \right] ,   \label{MILStrongDual}\\ 
& \quad \sum_{j \in \hat{\cJ}} y_j =1,   \label{MILRegularization}\\
& \quad \hat{\by} \ge \bzero,  \label{MILDualFeas2}\\
& \quad  \bx \in \mathbb{R}^n, \quad \bc \in \mathbb{R}^n,\quad \by \in \mathbb{R}^m_{\geq 0},\quad \hat{\bb} \in \mathbb{R}^{m_1} . \label{MILrange}
\end{align} 
\end{subequations}

In the above model, constraint \eqref{MILRegularization} enforces at least a subset of $\hat{\cJ}$ constraints to be binding for the learned solution $\bx$, meaning that the dual variables corresponding to the $\hat{\cJ}$ constraints for improvement can hold non-zero values, while are dual variables are forced to assume zero values. Other constraints of $\IL_g$ are similar to $\ILL$. Additionally, the second objective term in the $\IL_g$ formulation serves to choose a direction that provides the highest objective value (since $\bc' \bx = \bb' \by$ due to strong duality). For the settings where $\Omega$ is full-dimensional, non-empty, and free of redundant constraints, we first assert that $\IL_g$ is always feasible, then show that for the case where the observed decision $\bx_0$ is feasible for $\Omega$, the model has closed-form solutions for specific cases.

\begin{prop} \label{ILG_feasible}
      $(\hat{\ba}_i, \bx, \left [e_i; \,0 \right], \hat{\bb}_0)$ is feasible for $\IL_g$ where $\bx \in \Omega$ is on the boundary of $\Omega$ and binds constraint $i \in \hat{\cJ}$ (i.e. $\hat{\ba}_i \bx = \hat{b}^0_i$) and $e_i \in \mathbb{R}^{\left| \hat{\cJ} \right|}$ is the $i^{th}$ basis vector.
\end{prop}

Proposition \ref{ILG_feasible} indicates that $\IL_g$ has a feasible solution regardless of any condition on the observed decision $\bx_0$. Note that $\IL_g$ generalized $\ILL$ since if $\hat{\cJ} = \emptyset$, we arrive at the characterization of $\ILL$. For the case where $\hat{\cJ} \neq \emptyset$, however, $\IL_g$ provides an interpretable model to find realizations of $\FO$ that provide improved solutions in comparison to $\bx_0$. The next proposition characterizes the feasible set of $\IL_g$ and its relation to the feasible set of $\ILL$. We also discuss how the feasible set of $\IL_g$ changes as more constraints are added to $\hat{\cJ}$.

\begin{prop} \label{prop:ILG_feasiblity_properties}
    Let $\Phi (\hat{\cJ})$ be the feasible set of $\IL_g$ and let $\Phi$ be the feasible set of $\ILL$. Then we have the following:
    \begin{enumerate}
        \item $\Phi (\hat{\cJ})_{\bx} = \left\{ \bx \in \mathbb{R}^n | \bar{\bA} \bx \leq \bar{\bb}, \hat{\bb}_L\leq\hat{\bA} \bx\leq \hat{\bb}_U \right\}$ where $\Phi (\hat{\cJ})_{\bx}$ is the projection of $\Phi (\hat{\cJ})$ onto the space of $\bx$. 
        \item If $\hat{\cJ}_1$ and $\hat{\cJ}_2$ are two realizations of indices of desired constraints for improvement such that $\hat{\cJ}_2 \subseteq \hat{\cJ}_1$ , then $\Phi (\hat{\cJ}_2) \subseteq \Phi (\hat{\cJ}_1)$. Specifically, $\Phi \subseteq \Phi  (\hat{\cJ})$ for all $\hat{\cJ}$.
    \end{enumerate}
\end{prop}

The first part of Proposition \ref{prop:ILG_feasiblity_properties} discusses the properties of $\Phi (\hat{\cJ})$ by characterizing its projection on the space of the solutions. This results indicates the necessity of providing bounds on right-hand-side values in order to control the improvements. The second part discusses the relationship between the feasible sets of two realizations of $\IL_g$ with different sets of desired constraints for improvement, suggesting a trade-off in the number of right-hand-side elements subject to improvement and the level of control over the feasible set of the inverse problem. We next provide an important result of $\IL_g$, proving that improved solutions are obtained through $\IL_g$.



\begin{thm} \label{thm:ILG_solution}
    Let $(\bc, \bx, \by, \hat{\bb})$ be an optimal solution for $\IL_g({\bx_0, \bA, \bar{\bb}, \omega})$ and let $\bx_0 \in \Phi(\hat{\cJ})$. Then, $\bc' \bx \geq \bc'  \bx_0$.
\end{thm}

Theorem \ref{thm:ILG_solution} provides an important resultfor the general model $\IL_g$ in that regardless of the value of the objective weights $\omega$, utilization of $\IL_g$ provides a new solution that is guaranteed to perform at least as good as the observed solution. Additionally, the user can readily verify the level of improvement by comparing the right-hand-side values $\hat{\bb}$ with the original values $\hat{\bb}_0$ for each constraint in $\hat{\cJ}$. Also note that the only necessary condition for $\bx_0$ is feasibility for $\Phi(\hat{\cJ})$ rather than $\Phi$, showcasing model capability even for infeasible solution. This is particularly imporatant in our radiation therapy application since the constrained version of the problem might run into infeasibility issues for the observed plan. Finally, we note that in the case where $\bb$ is fully unknown, the problem has closed-form solutions based on the choice of objective weights, as shown in Proposition \ref{prop:ILg_Ab}.

\begin{prop} \label{prop:ILg_Ab}
Let $\bar{\cJ} = \emptyset$, then $\exists$ an optimal solution $\{\bc^*, \bx^*, \by^*, \hat{\bb}^* \} $ for $\IL_g({\bx_0, \bA, \bar{\bb}, 1})$ such that $\bx^* = \bx_0$ and $\exists$ an optimal solution $\{\bc^*, \bx^*, \by^*, \hat{\bb}^* \} $ $\IL_g({\bx_0, \bA, \bar{\bb}, 0})$ such that $\hat{\ba} \bx^* = \hat{b}$ for some $\hat{b} \in \{\hat{\bb}_L\} \cup \{\hat{\bb}_U\} $.
\end{prop}

Based on the results of Proposition \ref{prop:ILg_Ab}, for the case where all right-hand-side values are unknown, the model has extra levels of flexibility in providing realizations of $\FO$ that provide the desired effects. We note that, in general, such cases provide no additional information on whether the current decision $\bx_0$ is the best possible decision for a decision-making problem. Additionally, the general model $\IL_g$ has bilinear terms in the objective and in the constraints, making it hard to solve in applied settings. As such, in what follows, we only consider the case where exactly one of the right hand side parameters $\bb$ are uncertain or missing. For the case where $\left| \hat{\cJ} \right| = 1$, let $\hat{\bb} = \hat{b}$ and $\hat{\bA} = \hat{\ba}$. We simplify the optimization problem $\IL_g$ by setting $\hat{y}_j = 1$ for $ j \in \hat{\cJ}$ and $\bar{y}_j = 0$ $\forall j \in \bar{\cJ}$ such that $\sum_{j \in \cJ} y_j =1$. Employing the new values of the dual variables into $\IL_g$ yields a simpler form for the $\IL_g$ problem.
\begin{subequations} \label{MIL_detailed}
\begin{align}
\MIL({\bx_0, \bA, \bb,j}): \underset{\bx,\hat{b}}{\text{minimize}} & 
\quad \omega \cD(\bx, \bx_0) - (1-\omega) \hat{b}\\
\text{subject to} 
& \quad \bA \bx \leq \bb',   \label{MIOPrimal Feasblility1}\\ 
& \quad \ba^T \bx = \hat{b},  \label{MIL_Const_1} \\
& \quad  \bx, \bepsilon \in \mathbb{R}^n. \label{MIOrange}
\end{align}
\end{subequations}

\begin{figure}[htbp]
\begin{center}
\includegraphics[width =0.7 \linewidth]{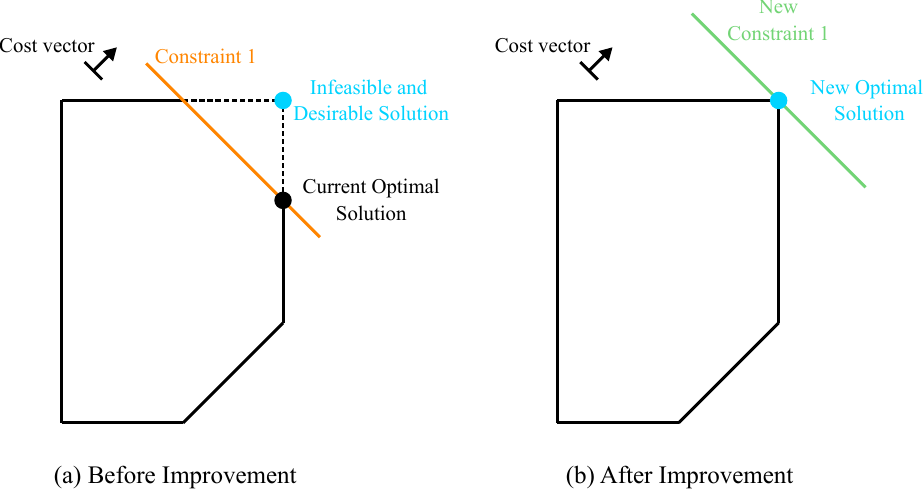}
\caption{\footnotesize A schematic example showcasing the improvement approach for existing solutions. (a) Before improvement, the current optimal solution is determined by the feasible region and the cost vector. An infeasible and desirable solution exists outside the feasible region. (b) After improvement, one of the constraints (Constraint 1) is allowed to move as much as possible in the direction of the known cost vector. This expands the feasible region, enabling the new optimal solution to achieve the desired improvements without compromising any of the remaining constraints.
} \label{Fig:methodology_example}
\end{center}
\end{figure}

Note that for formulation \eqref{MIL_detailed}, the second term of the objective function of $\ILL$ is changed to $(1-\omega) \hat{b}$, reducing this term to a linear objective. This is due to the fact that based on the strong duality condition, $\bc^T\bx = \hat{b}$. By defining $\hat{b}$ as a decision variable, the feasible set of $\FO$ is modified through the $\MIL$ model such that improved solutions are obtained in comparison to $\bx_0$. Noting the constraints of $\MIL$, the new learned solution is still feasible for all the remaining constraints due to condition \eqref{IOPrimal Feasblility1}. Additionally, constraint \eqref{MIL_Const_1} states that the new solution adheres to the new condition set forth by the modified value ($\hat{b}$). We note that this constraint is the strong duality constraint re-written using the new values of the parameters. While in our applications, $\bx_0$ consists of a single observed plan, $\MIL$ can be generalized to the case where multiple decisions over the same optimization problem are observed or multiple observations adhering to potentially similar patients can also be considered. The behavior of the multi-objective $\MIL$ model can be characterized by the relative values of the objective terms' weights. The main goal of MIL is to explore the potential for improving the right-hand-side terms $\bb$ and allowing for the first objective term to have larger weights makes the model not deviate from the observed, known previous decisions. This makes the model conscious about any desirable and unknown aspects of previous decisions to be maintained in the next iteration of predictions. Our assumption is that such an approach allows for increased acceptability of future plan recommendations for patients from such a model. Figure \ref{Fig:methodology_example} shows a schematic result of employing $\MIL$ to a simple example problem.

As shown in Figure \ref{Fig:methodology_example} demonstrates the improvement problem for the simple case of a two dimensional problem, in the case where only one constraint is considered for improving the observed decision, it is possible to assume that the cost vector is in the direction of the left hand side vector of the constraint and employ $\MIL$ to find potential improvement. We first note that as long as the original feasible set $\Omega (\bA, \bb)$ is non-empty, full dimensional and free of redundant constraints, $\MIL$ is feasible as it is a reduced version of $\IL_g$. Theorem \ref{thm:MIL_ILG} provides an important result with regard to optimal solutions of $\MIL$.

\begin{thm} \label{thm:MIL_ILG}
    If $\hat{\cJ} = \{j\}$ and $\{\bx^*, \hat{b}^* \} $ is optimal for $\MIL$, then there exists an optimal solution 
    $\{\bb, \bx, \by, \hat{\bb} \} $ for $\IL_g({\bx_0, \bA, \bar{\bb}, \omega})$ such that $\bx = \bx^*$ and $\hat{\bb} = \hat{b}^*$.
\end{thm}

Employing a similar method to different potential directions of improvement will result in a number of improved solutions that can be evaluated and chosen by an expert. Additionally, upon feasibility, multiple constraint improvement can also be considered.  In what follows, we show the results of applying the $\MIL$ formulation to a representative prostate cancer patient with a known, acceptable plan and show how $\MIL$ can provide improved objectives and plans for patients.

\section{ Radiation Therapy Treatment Plan Improvement Methodology} \label{sec:RT_methodology}
Based on the decision improvement scheme we developed in \cref{sec:IL_methodology}, in this section, we aim to employ the methods to develop an optimization model that is capable of showcasing how much improvement is possible for a certain objective in the radiation therapy treatment planning problem. We first discuss the treatment planning optimization problem and characterize a linear version of the radiation therapy treatment planning problem. We employ known techniques to model dose-volume objectives as linear constraints to the model as well. We then apply the $\MIL$ model to the treatment planning problem where we assume a plan is given in addition to a ``direction of improvement''.

\subsection{Treatment Plan Optimization and Dose-Volume Constraints}
In this section, we detail the treatment plan optimization problem where an optimal beam intensity profile is desired to render an optimal dose-distribution.
We consider the setting where a dose distribution is desired based on available knowledge of the three dimensional structure of the patient and a set of clinical objectives for the target and non-target structures. In this setting, let $\bD$ be the matrix of doses corresponding to different voxels in the three dimensional body. Let $w_b$ represent the intensity of beam $b$ from the set of all beams $\mathcal{B}$ and the set of all intensities $\bW$ and let $D_{i,b}$ denote the dose-influence matrix that represents the dose delivered to voxel $i$ from beam $b$ at unit intensity. We borrow some of the notations in defining the treatment planning optimization problem from previous work \citep{brady2006new}, however, we specifically consider dose-volume constraints in the definition of the radiation therapy treatment planning problem. As such, in general, the radiation therapy treatment planning optimization problem can be formulated as follows.
\begin{subequations} \label{formulation:IPP}
\begin{align}
\IPP: 
\underset{\bD, \bW}{\text{minimize}} & \quad \sum_{s \in \cS} w_{s} \bF(s)
\label{IPP_obj}\\
\text{subject to} 
& \quad \bF(s) = \sum_{v \in s} \cD (d_v, d_P^s),   \quad \forall s \in \cS \label{const:IPP_dose_deviation}\\
& \quad \max_{v \in s_{\alpha}} d_v \leq U_s^{\alpha},   \quad \forall \alpha \in \bA_s, \forall s_{\alpha} \subseteq s \ where \ \left|s_{\alpha} \right| = \alpha \left|s \right|, \forall s \in \cS \label{const:IPP_DV_UB}\\
& \quad \min_{v \in s_{\beta}} d_v \geq L_s^{\beta},   \quad \forall \beta \in \bB_s, \forall s_{\beta} \subseteq s \ where \ \left|s_{\beta} \right| = \beta \left|s \right|, \forall s \in \cS \label{const:IPP_DV_LB}\\
& \quad d_v = \sum_{b \in \mathcal{B}} D_{v,b} \cdot w_b, \quad \forall v \in s, \quad \forall s \in S \label{eq:rtp_dose_calc} \\
& \quad w_b \geq 0, \quad \forall b \in \mathcal{B}\ \label{eq:rtp_nonneg_beam}\\
& \quad d_v \geq 0,   \quad  \forall v \in s,\forall s \in \cS 
\label{const:IPP_non-negativity}
\end{align}
\end{subequations}

In the above model, $\bF$ in constraint \eqref{const:IPP_dose_deviation} represents the objective functions for minimizing deviations of voxel doses to the prescribed dose ($d_P^s$) for structure $s$. The main objective of $\IPP$ is then to minimize a weighted sum of all such deviation functions for different structures in the problem. Equally important, constraints \eqref{const:IPP_DV_UB} and \eqref{const:IPP_DV_LB} represent dose-volume criteria specified for given sets of volume fractions. These sets are denoted as $\bA_s$ and $\bB_s$ for constraints \eqref{const:IPP_DV_UB} and \eqref{const:IPP_DV_LB} respectively. Finally, constraint \eqref{eq:rtp_nonneg_beam} connects the beam intensities to voxel doses through a known dose-influence matrix and \eqref{const:IPP_non-negativity} asserts that dose levels for all voxels have to be non-negative. The radiation therapy treatment planning problem includes inherent trade-offs in the objective, mainly on exposing the target structure to higher levels of radiation and sparing the organs at risk (e.g. rectum, bladder and penile bulb for the case of prostate cancer). Assuming equivalent soft constraints in terms of objective terms, for any realization of the parameters $\bA_s$ and $\bB_s$ along with bounds $U$ and $L$, there can be multiple solutions that provide desirable effects due to the nature of the problem. We show in subsequent sections that our plan improvement methodology allows for exploring previously non-attainable solutions by updating the parameters of $\IPP$ using existing plans. Additionally, note that $\IPP$ is a multi-objective optimization problem in general.

A linear programming model for the radiation therapy treatment planning problem ($\IPP$) can be found in  \cite{romeijn2003novel}. In this work,  we consider a linear version of $\IPP$ that is capable of including linear approximations of dose-volume criteria for both the target and non-target structures in the form of conditional value at risk constraints, if such criteria are specified. Formulation \ref{formulation:IPP_detailed} provides a detailed account of $\IPP$ as a linear optimization problem including dose-volume criteria using notions from previous literature \citep{mahmoudzadeh2015robust}.
\begin{subequations} \label{formulation:IPP_detailed}
\begin{align}
\underset{d, u, o, z, m, w}{\text{minimize}} & \quad \sum_{N \in \mathcal{O}} (\sum_{i \in N} (\delta_i u_i^N + \zeta_i o_i^N)) + \sum_{N \in \mathcal{O}} ( \eta_N z^N +  \kappa_N m^N) + \sum_{T \in \mathcal{T}} (\sum_{i \in T} (\lambda_i u_i^T + \mu_i  o_i^T))
\label{eq:IPP_obj}\\
\text{subject to} 
& \quad u_i^T \geq d_i^{PT} -d_i^{T},   \quad \forall i \in T,  \forall T \in \mathcal{T} \label{constr_underdose_target}\\
& \quad o_i^T \geq d_i^{T} -d_i^{PT},   \quad \forall i \in T,  \forall T \in \mathcal{T} \label{constr_overdose_target}\\
& \quad u_i^N \geq d_i^{PN} -d_i^{N},   \quad \forall i \in N,  \forall N \in \mathcal{O} \label{constr_underdose_non_target}\\
& \quad o_i^N \geq d_i^{N} -d_i^{PN},   \quad \forall i \in N,  \forall N \in \mathcal{O} \label{constr_overdose_non_target}\\
& \quad z_N \geq d_i^N,   \quad \forall i \in N,  \forall N \in \mathcal{O} \label{constr_maxdose_non_target}\\
& \quad m_n = \frac{1}{\left | N \right |} \sum_{i \in N}d_i,   \quad  \forall N \in \mathcal{O}
\label{constr_meandose_non_target}\\
& \quad  \gamma_{\alpha}^N + \frac{1}{(1-\alpha) \left | N \right |}  \sum_{i \in N} \Bar{d_i}^{N,\alpha} \leq  U^{\alpha}_N,   \quad  \forall \alpha \in \mathcal{U}^N  \label{IPP_DVH_const_1} \\
& \quad  \gamma_{\alpha}^N + \frac{1}{(1-\alpha) \left | N \right |}  \sum_{i \in N} \Bar{d_i}^{N,\alpha} \geq  L^{\alpha}_N,   \quad  \forall \alpha \in \mathcal{L}^N  \label{IPP_DVH_const_3} \\
& \quad \Bar{d_i}^{N,\alpha} \geq  d_i - \gamma_{\alpha}^N,   \quad \forall i \in N,  \forall N \in \mathcal{O} \label{IPP_DVH_const_2} \\ 
& \quad d_i = \sum_{b \in \mathcal{B}} D_{i,b} \cdot w_b,   \quad  \forall i \in S \label{dose_influence_constraint}\\
& \quad \Bar{d_i} \geq 0,   \quad  \forall i \in S \\
& \quad w_b \geq 0, \quad \forall b \in \mathcal{B}\ \label{eq:IO_nonneg_beam}\\
& \quad d_i \geq 0,   \quad  \forall i \in S 
\label{non-negativity}
\end{align}
\end{subequations}

In formulation \ref{formulation:IPP_detailed}, constraints \eqref{constr_underdose_target} and \eqref{constr_overdose_target} measure the overdose and underdose to the target structures with respect to the prescribed dose levels, constraints \eqref{constr_underdose_non_target} and \eqref{constr_overdose_non_target} measure the overdose and underdose to the non-target structures with respect to the prescribed dose levels, constraint \eqref{constr_maxdose_non_target} limits the maximum dose levels to non-target structures, constraint \eqref{constr_meandose_non_target} measures the average dose to non-target structures. Constraints \eqref{IPP_DVH_const_1} and \eqref{IPP_DVH_const_2} represent the conditional value at risk linear approximations for dose volume constraints \citep{mahmoudzadeh2015robust} and constraint \eqref{dose_influence_constraint} represents the dose influence constraint while constraints \eqref{non-negativity} ensure non-negativity of voxel doses. Therefore, the objective of formulation \ref{formulation:IPP_detailed} is a weighted sum of the respective measured terms in the constraints. \cref{table:IPP_descriptions} provides a detailed account of the terms used in formulation \ref{formulation:IPP_detailed}.

\begin{table}[htbp]
\renewcommand{\arraystretch}{1.5}
\caption{Descriptions of the terms used in the definition of the linear radiation therapy treatment planning optimization problem in formulation \ref{formulation:IPP_detailed}.} \label{table:IPP_descriptions}
\begin{tabular}{>{\centering}m{0.1\textwidth}>{\arraybackslash}m{0.80\textwidth}}
\hline
\textbf{Term} &  \textbf{Description}  \\ 
\hline 
$u_i^T$ &  Underdose to voxel $i$ of the target structure $T$  \\ 
$o_i^T$ &  Overdose to voxel $i$ of the target structure $T$  \\ 
$u_i^N$ &  Underdose to voxel $i$ of the    \\ 
$o_i^N$ &  Overdose to voxel $i$ of the non-target structure $N$  \\ 
$z^N$ &  Maximum voxel dose to the non-target structure $N$  \\ 
$m^N$ &  Average dose to the non-target structure $N$ \\ 
$d_i^{PT}$ &  Prescribed dose for voxel i of the target structure $T$  \\ 
$d_i^{T}$ &  variable of the dose for voxel i of the target structure $T$  \\ 
$d_i^{PN}$ &  Prescribed dose for voxel i of the non-target structure $N$   \\ 
$d_i^{N}$ &  variable of the dose for voxel i of the non-target structure $N$  \\ 
$\mathcal{T}$ &  Set of all target structures  \\ 
$\mathcal{O}$ &  Set of all non-target structures  \\ 
$\gamma_{\alpha}^N$, $\Bar{d_i}^{N,\alpha}$  &  Auxiliary variables for the DVH constraints  \\ 
$U^{\alpha}_N$ &  Upper dose limit of the DVH constraint for fractional volume level $\alpha$ in the non-target structure $N$ \\ 
$L^{\alpha}_N$ &  Lower dose limit of the DVH constraint for fractional volume level $\alpha$ in the non-target structure $N$  \\ 
$\mathcal{U}^N$ &  Set of fractional volumes $\alpha$ for which an upper limit DVH constraint is defined  \\ 
$\mathcal{L}^N$ &  Set of fractional volumes $\alpha$ for which an lower limit DVH constraint is defined  \\ 
$d_i$ &  Dose variable to voxel $i$  \\ 
$\mathcal{B}$  &  Set of all beams  \\ 
$D_{i,b}$ &  Dose-fluence matrix parameter  \\ 
$w_b$ &  Variable for the intensity of beam $b$  \\ 
\hline
\end{tabular}
\end{table}

In what follows, we consider the linear formulation \ref{formulation:IPP_detailed} as the treatment plan optimization problem. We next aim to detail the key methodology block that we employ to improve existing radiation therapy treatment plans by proposing new right hand side values for some of the bounds for the dose-volume constraints.

\subsection{Treatment Plan Improvement Methodology}
\begin{figure}[htbp]
\begin{center}
\includegraphics[width =0.7 \linewidth]{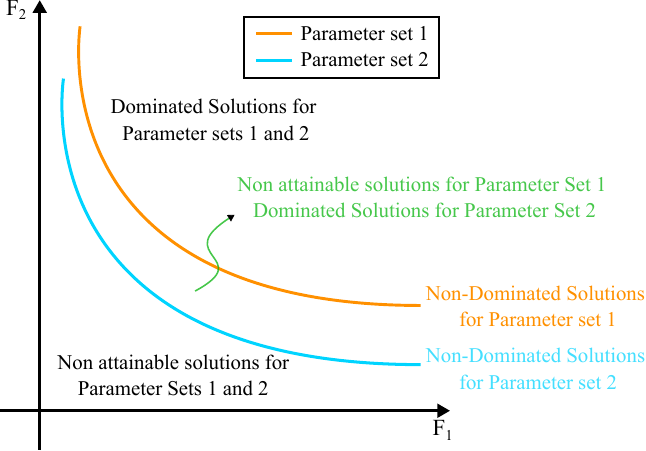}
\caption{\footnotesize Schematic representation of the Pareto frontier for a multi-objective optimization problem with competing objectives F1 and F2. The figure illustrates how the Pareto frontier changes for different realizations of the problem parameters. Parameter set 1 and Parameter set 2 represent two distinct sets of parameters that define the feasible region of the optimization problem. The Pareto frontier for Parameter set 1 differs from the Pareto frontier for Parameter set 2 , indicating that the choice of parameters affects the trade-offs between the objectives and the attainable solutions. Solutions that are dominated by other solutions for a given parameter set are represented, while non-attainable solutions lie outside the feasible region for both parameter sets.} \label{Fig:methodology_Pareto}
\end{center}
\end{figure}

Using the methodology developed in \cref{sec:IL_methodology} as a core approach for automatizing the treatment plan modification and improvement, one can update the treatment planning framework in \cref{Fig:IPP_framework} to realize an automated treatment process for iteratively improving treatment plan. For the case of radiation therapy treatment planning, dose volume constraints are widely used in forming the optimization problem. For such constraints that can be formulated as shown in constraints \eqref{IPP_DVH_const_1}, \eqref{IPP_DVH_const_2}, and \eqref{IPP_DVH_const_3}, for different regions of interest, the direction of improvement (cost function) is usually known to the planner. For example, for a dose limit constraint for a certain organ at risk, smaller left-hand-side values are desirable. As such, reducing the dose limit such that feasible solutions exist that adhere to the new limit is desirable and provides improved solutions. \cref{Fig:methodology_Pareto} provides a schematic figure of the effect of improving solutions. As shown, by providing modified realizations of some of the right hand side parameters of the constraints, new improved solutions can be achieved.

Considering formulation \eqref{formulation:IPP_detailed}, in order to explore potential plan improvements, we assume a single plan is given alongside with the objectives or constraints that generate the plan. Applying the $\MIL$ formulation to the linear version of the radiation therapy treatment planning problem provided in formulation \eqref{formulation:IPP_detailed} yields the following formulation where an improved upper-bound for a selected dose volume constraint for a specific organ at risk is achieved. 
\begin{subequations} \label{formulation:MIL_IPP_detailed}
\begin{align}
\underset{d, u, o, z, m, w, \bE, U_{\hat{N}}}{\text{minimize}} & \quad \omega \cD(\bE) - (1-\omega) U_{\hat{N}}
\label{eq:MIL_IPP_obj}\\
\text{subject to} 
& \quad  \gamma_{\alpha}^N + \frac{1}{(1-\alpha) \left | N \right |}  \sum_{i \in N} \Bar{d_i}^{N,\alpha} =  U_{\hat{N}},     \label{MIL_IPP_DVH_const_1} \\
& \begin{bmatrix} d \\ u \\ o \\ z \\ m \\ w \\ \end{bmatrix} = 
\begin{bmatrix} d_0 \\ u_0 \\ o_0 \\ z_0 \\ m_0 \\ w_0 \\ \end{bmatrix} -
\begin{bmatrix} \bepsilon_d \\ \bepsilon_u \\ \bepsilon_o \\ \bepsilon_z \\ \bepsilon_m \\ \bepsilon_w \\ \end{bmatrix} , \label{MIL_IPP_IOOnePoint_1} \\
& \quad \eqref{constr_underdose_target} - \eqref{non-negativity}.
\end{align}
\end{subequations}

Formulation \eqref{formulation:MIL_IPP_detailed} is the rendition of $\MIL$ for the $\IPP$ problem. In this formulation, a new constraint is added to formulation \eqref{formulation:IPP_detailed} where a certain dose-volume constraint is to be improved for the organ at risk ${\hat{N}}$. Other constraints are similar to the previous formulation. Using formulation \eqref{formulation:MIL_IPP_detailed} iteratively, one can uncover potential improvement for an observed treatment plan automatically. The inverse optimization model $\MIL$ is employed in a loop to constantly improve the given treatment plans until a certain criterion is met.

After an initial (current) plan is achieved using a treatment planning optimization model through existing patient data and a set of initial objectives, including DV objectives, the inverse optimization model is used to provide an improvement to the current plan, resulting in an improved plan. If the improved plan varies insignificantly from the current plan based on some user-specified criterion, the process can be repeated until possible mathematical improvements to the plan are uncovered. As previously mentioned, these improvements can include better sparing of some organs at risk or even better/uniform coverage of the target. In what follows, we employ a retrospective analysis on previously treated patients to provide use cases of this methodology for improving existing treatment plans. Naturally, the same procedures can be used for future patients, given a set of initial objectives that arise from expert knowledge or data-driven models.

\section{Results} \label{sec:results}

We provide summary results of applying the $\MIL$ model to four prostate cancer patients. We first detail the data required for the models in \cref{sec:results_data}. We then provide results for applying the $\MIL$ model for improvement of a single OAR's dose volume curve by improving a single constraint's dose limit in \cref{sec:results_signluar_OAR} and finally touch on implications on the pareto frontier of the multi-objective representation of the treatment planning problem before and after improving a single dose limit in \cref{sec:results_Pareto}.

\subsection{Data} \label{sec:results_data}
The proposed model is applied to four clinical prostate cancer patient cases to evaluate its ability to improve upon existing clinically accepted IMRT treatment plans. The patient cohort consisted of two standard prostate cancer cases (Patients 1 and 2), where the planning target volume (PTV) was largely confined within the prostate gland, as well as two cases requiring craniocaudal cone-down plans (Patients 3 and 4) due to the proximity of the seminal vesicles to the rectum. For each patient, the clinically delivered dose distribution and the corresponding predicted objective list that replicated the dose distribution were available. \cref{table:rois} indicates the regions of interest for each patient that are considered in the optimization models $\MIL$. The data is sourced from DICOM files, including Computed Tomography (CT) images and contouring information. In addition to these inputs, the models need an existing plan for each patient alongside information on the objective list that generates the existing plan. Considering that such information is not unified and sometimes not available, we used a machine learning-based objective prediction model to replicate the plan \citep{wu2011data}. \cref{Table:Twain_Plan_Objectives} shows a list of the predicted objectives for a representative patient. Our goal is to evaluate possible plan by improving the objective right-hand-side values.

\begin{table}[]
\renewcommand{\arraystretch}{2}
\caption{List of ROIs considered in the optimization for each patient. } \label{table:rois}
\begin{center}
\begin{tabular}{c|c|c|c|c}

\textbf{}     & \textbf{Patient 1} & \textbf{Patient 2} & \textbf{Patient 3} & \textbf{Patient 4} \\ \hline \hline
\textbf{ROIs} & \makecell{PTV \\ Penile Bulb\\Bladder\\Rectum\\Sigmoid\\Left Femoral Head\\Right Femoral Head}   & \makecell{PTV\\Bladder\\Rectum\\Left Femoral Head\\Right Femoral Head\\}   & \makecell{PTV\\PTVNodes\\SmallBowel\\Bladder\\Rectum\\Sigmoid\\Left Femoral Head\\Right Femoral Head}  & \makecell{PTV\\Bladder\\Penile Bulb\\Rectum\\Sigmoid\\Left Femoral Head\\Right Femoral Head\\} \\
\hline\hline
\end{tabular}
\end{center}
\end{table}

\begin{table}[htbp]
\begin{center}
\caption{Predicted dose-volume objectives and their associated parameters used to generate the clinically acceptable treatment plan for the representative prostate cancer patient.}
\label{Table:Twain_Plan_Objectives}
\begin{tabular}{|c|c|c|c|c|c|}
\hline
\textbf{Objective type} & \textbf{ROI Name}                & \textbf{ROI Type} & \textbf{Dose} & \textbf{Weight} & \textbf{Percentage} \\
\hline \hline
Min DVH                 & PTV                              & Target            & 3625             & 100              & 100              \\
Max Dose                & PTV                              & Target            & 3700             & 100              & \_               \\
Uniform Dose            & PTV                              & Target            & 3625             & 200              & \_                 \\
Max DVH                 & Penile Bulb                      & OAR               & 607              & 1                & 99               \\
Max DVH                 & Bladder                          & OAR               & 3625             & 1                & 100              \\
Max DVH                 & Rectum                           & OAR               & 3625             & 5                & 100              \\
Max DVH                 & Sigmoid                          & OAR               & 79               & 1                & 99               \\
Max DVH                 & Left Femoral Head                   & OAR               & 1174             & 1                & 100              \\
Max DVH                 & Right Femoral Head                   & OAR               & 1134             & 1                & 100              \\
Max DVH                 & Left Femoral Head                     & OAR               & 1030             & 1                & 95               \\
Max DVH                 & Right Femoral Head                   & OAR               & 991              & 1                & 95               \\
Max DVH                 & Bladder                          & OAR               & 246              & 1                & 50               \\
Max DVH                 & Bladder                          & OAR               & 1300             & 1                & 80               \\
Max DVH                 & Bladder                          & OAR               & 2341             & 1                & 90               \\
Max DVH                 & Bladder                          & OAR               & 3009             & 1                & 95               \\
Max DVH                 & Rectum                           & OAR               & 244              & 1                & 50               \\
Max DVH                 & Rectum                           & OAR               & 1345             & 5                & 80               \\
Max DVH                 & Rectum                           & OAR               & 2541             & 5                & 90               \\
Max DVH                 & Rectum                           & OAR               & 3163             & 5                & 95               \\
Max DVH                 & Penile Bulb                      & OAR               & 571              & 1                & 97               \\
Max Dose                & Penile Bulb                      & OAR               & 626              & 1                & \_                 \\
Max DVH                 & Rectum\_ring\_0-10mm         & OAR               & 3625             & 1                & 100              \\
Max DVH                 & ptv\_ring\_10-20mm & OAR               & 2258             & 1                & 100              \\
Max DVH                 & ptv\_ring\_20-40mm & OAR               & 2045             & 1                & 100              \\
Max DVH                 & ptv\_ring\_40-60mm & OAR               & 1823             & 1                & 100              \\
Max DVH                 & ptv\_ring\_0-10mm  & OAR               & 2767             & 1                & 100              \\
Max DVH                 & ptv\_ring\_0-10mm  & OAR               & 2261             & 1                & 80               \\
Max DVH                 & ptv\_ring\_10-20mm & OAR               & 1529             & 1                & 80               \\
Max DVH                 & ptv\_ring\_20-40mm & OAR               & 1219             & 1                & 80               \\
Max DVH                 & ptv\_ring\_40-60mm & OAR               & 952              & 1                & 80               
\\
\hline           
\end{tabular}
\end{center}
\end{table}

\subsection{Simulating Plans with matRad}
The open source software matRad \cite{wieser2017development} is used for intensity modulated radiation therapy (IMRT) treatment planning. In order to provide fair comparison between the clinical plan and the improved plan generated from the proposed models in this work, we use matRad to first replicate the clinical plan using the predicted objectives and a fixed beam configuration across the board. We then use the resulting dose distribution as the available plan for the proposed models, which are also rendered in matRad. In order to construct $\IPP$, we use the known objectives used to arrive at the available plan using matRad on a multi-objective optimization scheme. For all plans, we assume only one fraction for all treatment plans for simplicity and consider nine equally distanced fixed beams. All the analyses are done using a conventional computer equipped with a 2 GHz Quad-Core Intel Core i5 processor with 16 GB 3733 MHz LPDDR4X memory using the latest version of Gurobi (9.5) \citep{gurobi}.

\subsection{Improvement on Singular OAR} \label{sec:results_signluar_OAR}

We primarily focus on improving the rectum dose-volume curve by updating the value of the dose limit for the 30\% fractional volume using the inverse process $\MIL$. To test generalization to other OARs, we replicate the analysis for the case of other important OARs for one of the patients as well. We explore improving dose levels for a single OAR using the available plan's dose distribution and known objectives for all patients.  For all cases, we provide improvements over the 30\% fractional volume level of the OAR that is subject to improvement. If a dose-volume objective already exists for this criteria, its dose limit is treated as an unknown parameter. Otherwise, a new dose-volume constraint is added with an unknown dose limit. 

\begin{figure}[htbp]
\includegraphics[width =0.9 \linewidth]{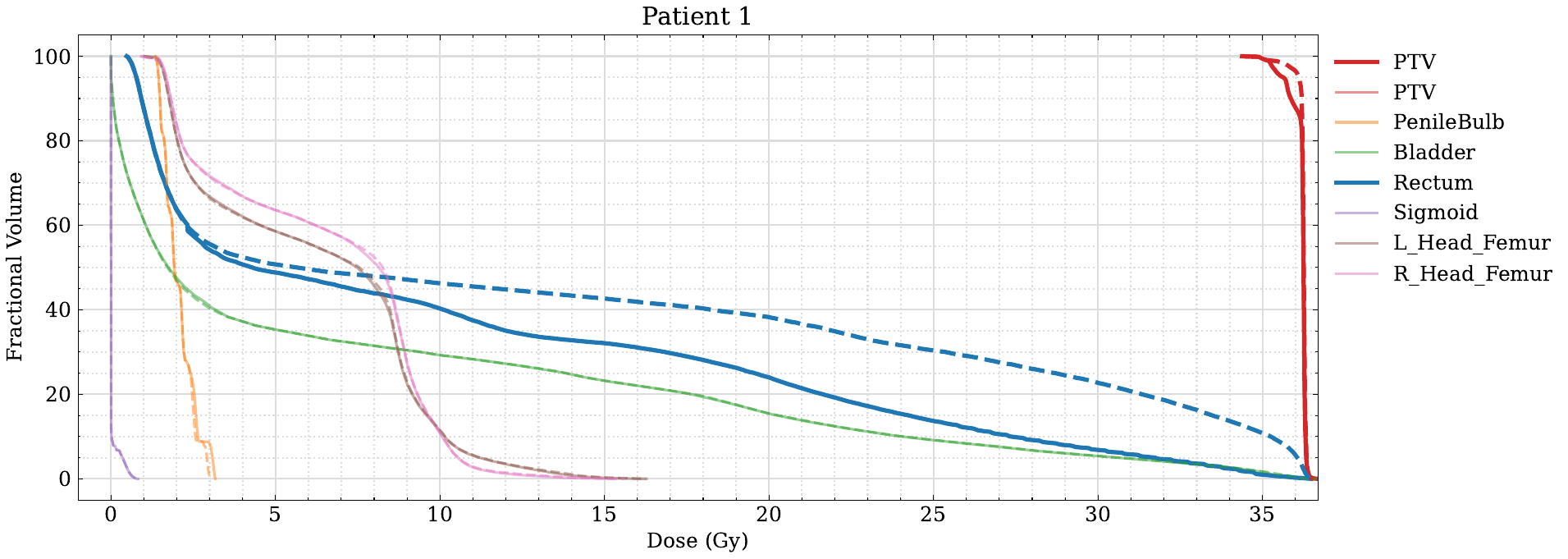}
\includegraphics[width =0.98 \linewidth]{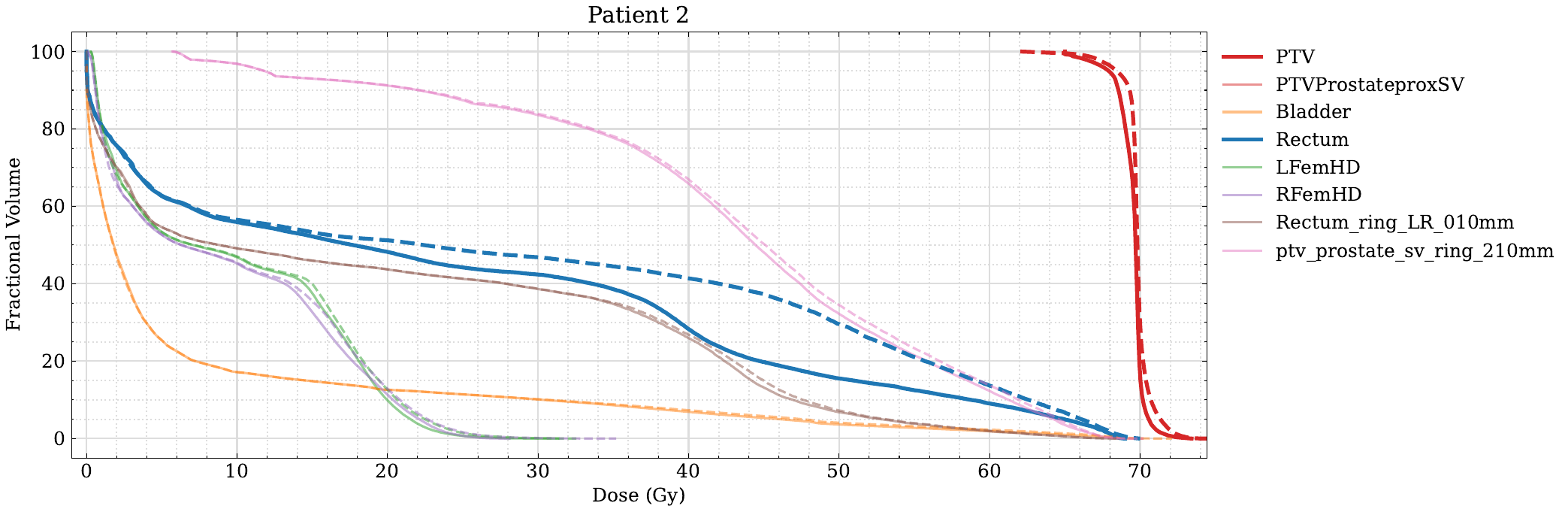}
\includegraphics[width =0.91 \linewidth]{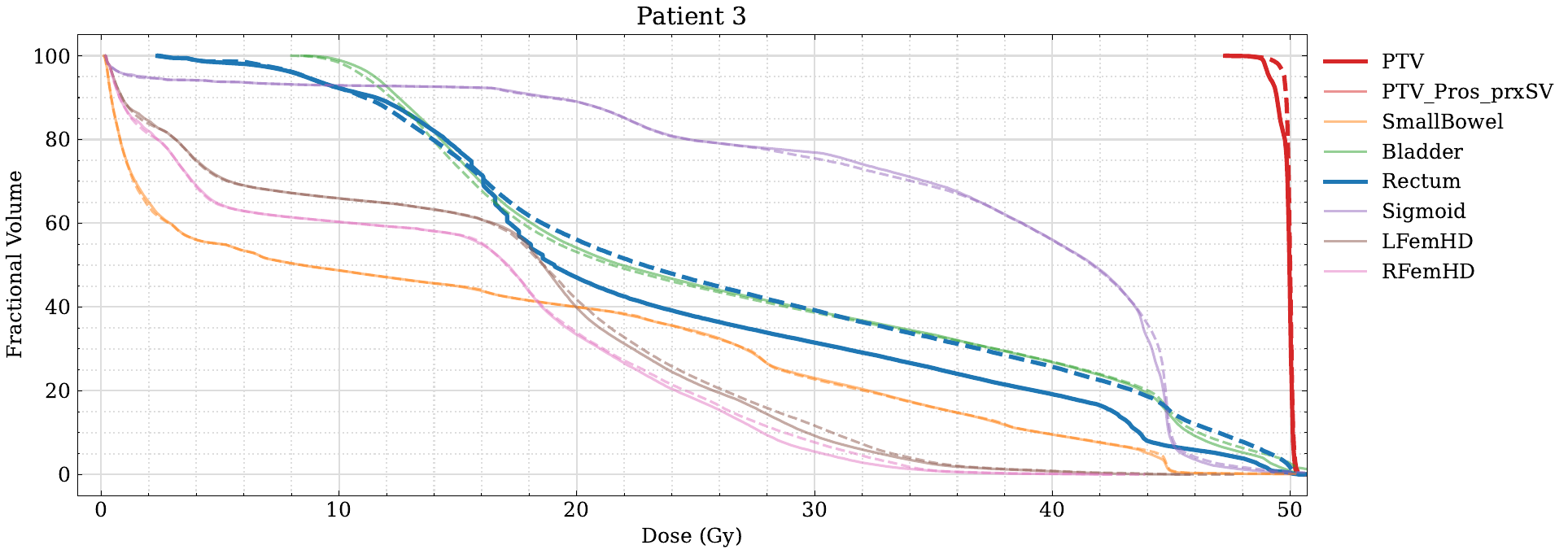}
\includegraphics[width =0.9 \linewidth]{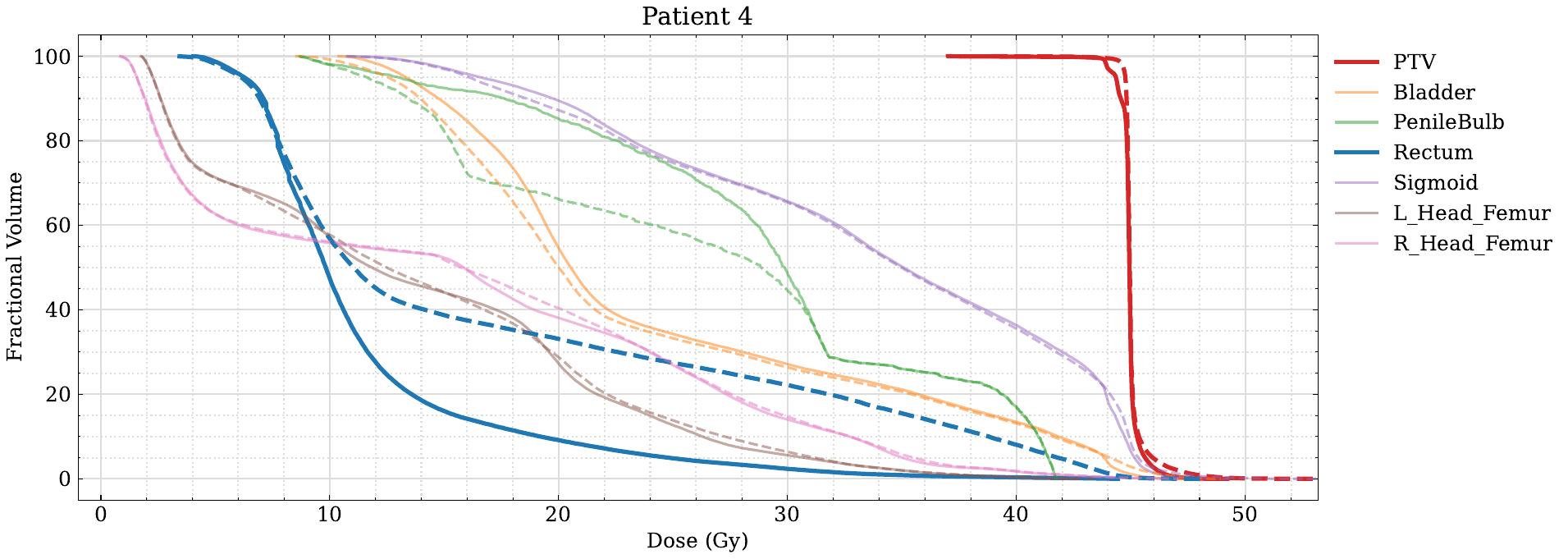}
\caption{\footnotesize DVH curves for all patients before (dashed lines) and after (solid lines) improvement (optimization results). As shown, after running the $\MIL$ model for improving rectum dose levels, dose reductions are seen for all patients on the OAR in question while target dose levels remain roughly unchanged.} \label{Fig:DVH_comparison_rectum}
\end{figure}

\cref{Fig:DVH_comparison_rectum} compares the DVH curves of the original and $\MIL$-improved plans for all four patients. Improvements are seen across the board for all patients by reducing the overall dose levels to rectum while keeping the dose levels to the target roughly unchanged. For each patient, the DVH curves demonstrate a notable reduction in rectum dose levels in the MIL-optimized plans compared to the original plans. The rectum DVH curves in the MIL plans consistently shift towards lower dose values, indicating improved rectum sparing across all four cases. This improvement is particularly evident in the intermediate to high dose regions, where the MIL plans achieve substantial reductions in the percentage of rectum volume receiving higher doses.

Importantly, the DVH curves for the PTVs remain nearly identical between the original and MIL-optimized plans for all patients. The close overlap of the PTV curves indicates that the MIL model successfully maintains target coverage while improving rectum sparing. This finding is crucial, as it demonstrates that the dosimetric improvements in the rectum are not achieved at the expense of compromising the primary treatment objectives. In most cases, the DVH curves for OARs such as the bladder, penile bulb, and femoral heads exhibit minimal changes between the original and MIL plans. This consistency suggests that the MIL model effectively focuses on improving rectum sparing without introducing significant dosimetric trade-offs in other critical structures. While the magnitude of rectum dose reduction varies among patients, the consistent improvement in rectum DVH curves across all cases highlights the robustness and adaptability of the MIL approach. Although the DVH curves for other OARs generally showed minimal changes, it is worth noting that for Patient 4, the penile bulb dose increased in the MIL-optimized plan compared to the original plan. However, our clinical collaborator suggested that the substantial improvement in rectum sparing achieved by the MIL plan outweighs the potential drawbacks of the increased penile bulb dose. The collaborator emphasized that the overall benefit of the MIL-optimized plan, particularly in terms of reducing the risk of radiation-induced rectal toxicities, justifies its consideration for clinical implementation.

\begin{figure}[htbp]
\begin{center}
\includegraphics[width =0.3 \linewidth]{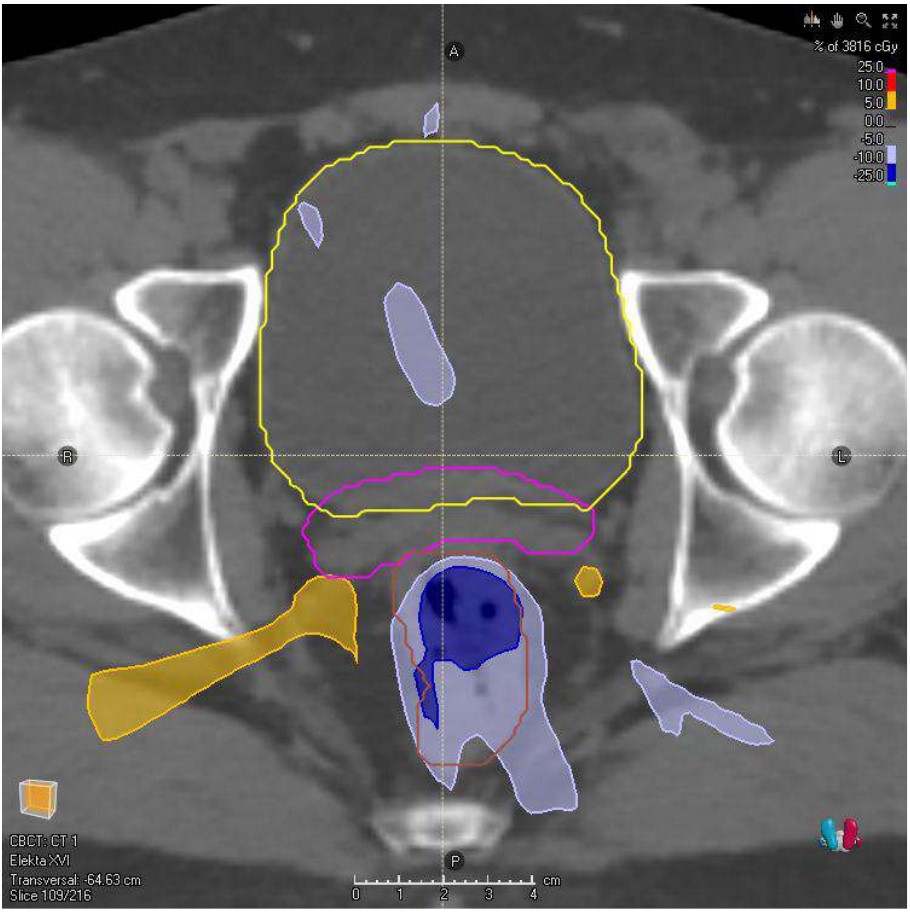}
\includegraphics[width =0.3 \linewidth]{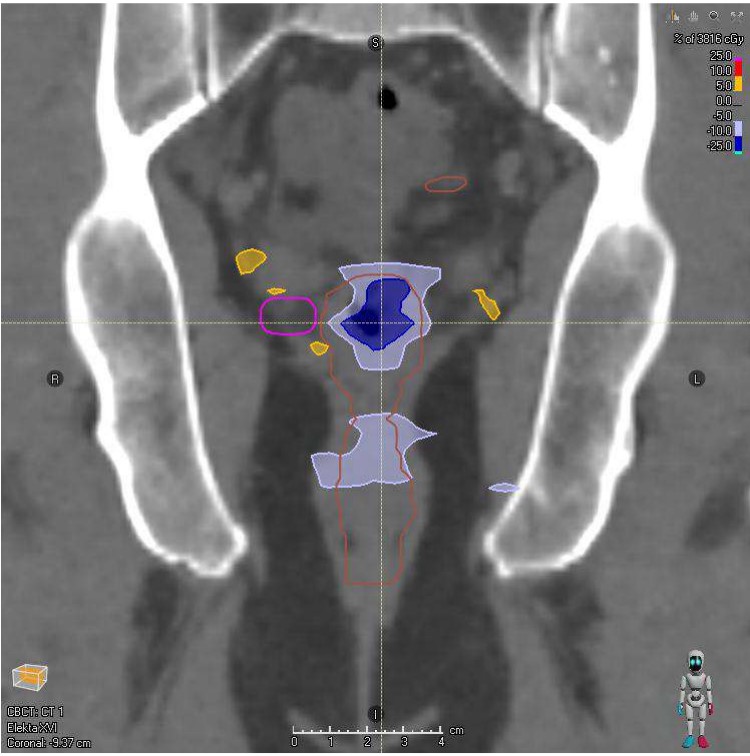}
\includegraphics[width =0.3 \linewidth]{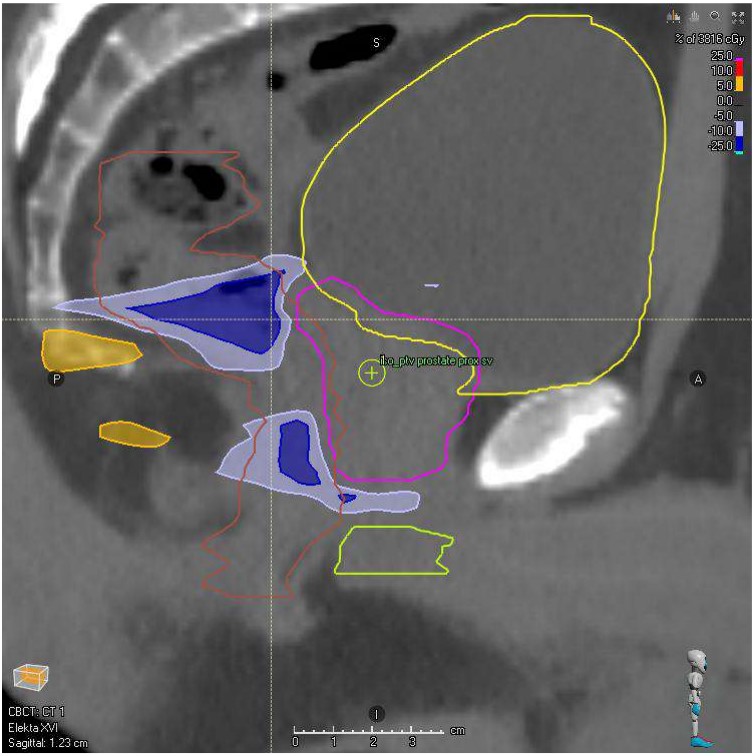}
\caption{\footnotesize Spatial dose distribution comparisons between the original (left) and MIL-optimized (right) plans for each patient, showcasing dose differences in the axial, sagittal, and coronal planes. The MIL plans exhibit notable reductions in rectum dose levels while maintaining similar dose distributions for the planning target volume  and other organs-at-risk. These improvements are achieved through subtle modifications in the dose profiles guided by the updated optimization parameters, demonstrating the effectiveness of the MIL approach in enhancing rectum sparing without compromising target coverage or overall plan quality.} \label{Fig:dose_profile_comparison}
\end{center}
\end{figure}

\cref{Fig:dose_profile_comparison}  showcases the spatial dose distribution differences between the original treatment plans and the improved plans generated by the $\MIL$ model for a representative patient. The figure consists of three dose distribution cuts for the patient, highlighting the nuances in dose profiles that lead to improvements in the rectum dose levels due to changes in the optimization parameters, even with fixed beam configurations. The dose distribution cuts reveal that the $\MIL$-optimized plans exhibit notable differences compared to the original plans. These differences are particularly evident in the regions surrounding the rectum, where the $\MIL$ plans demonstrate reduced dose levels. Importantly, the dose distributions for the target volumes  remain largely unchanged between the original and $\MIL$-optimized plans. This consistency in PTV dose coverage indicates that the MIL model successfully maintains the primary treatment objectives while improving rectum sparing. The dose distribution cuts also reveal that the improvements in rectum sparing are not associated with substantial changes in dose levels to other nearby OARs. 

\cref{table:rectum_comparison} shows a more detailed account of the differences between the improved plan using $\MIL$ and the observed plan. As shown in the table, the reduction in rectum dose-levels are around 20 to 30 percent while changes in target coverage remain fairly insignificant. Use of the improvement methodology additionally results in new recommendations for dose-level values, resulting in improved objectives that result in objectively better plans. As such, the  shows the values of the parameter subject to improvement using $\MIL$ for different OARs. 

\begin{table}[htbp]
\begin{center}
\renewcommand{\arraystretch}{1.5}
\caption{Comparison of key dose-volume metrics between the MIL-optimized (proposed) and original (observed) treatment plans for the rectum and planning target volum across all four prostate cancer patients.} \label{table:rectum_comparison}
\begin{tabular}{c|cc|cc|cc|cc}
\multirow{3}{*}{} & \multicolumn{2}{c|}{Patient 1}                                   & \multicolumn{2}{c|}{Patient 2}                                   & \multicolumn{2}{c}{Patient 3}                                   & \multicolumn{2}{c}{Patient 4}                                    \\ \cline{2-9} 
                  & $\MIL$    & \multicolumn{1}{c|}{Observed }          & $\MIL$    & \multicolumn{1}{c|}{Observed }           & $\MIL$    & \multicolumn{1}{c|}{Observed } & $\MIL$    & \multicolumn{1}{c}{Observed }         \\ \hline \hline
\textbf{$D^{PTV}_{max}$}              & 36.68 & \multicolumn{1}{c|}{36.68}         & 80.11 & \multicolumn{1}{c|}{81.11}         & 73.41 & \multicolumn{1}{c|}{74.43}     & 49.21 & \multicolumn{1}{c}{53.17}     \\
\textbf{$D^{PTV}_{95\%}$}               & 36.32 & \multicolumn{1}{c|}{36.33}         & 80.64 & \multicolumn{1}{c|}{80.7}          & 70.56 & \multicolumn{1}{c|}{71.33}       & 45.61 & \multicolumn{1}{c}{45.94} \\
\textbf{$D^{PTV}_{mean}$}              & 36.06 & \multicolumn{1}{c|}{36.22}         & 79.95 & \multicolumn{1}{c|}{80.2}         & 69.49 & \multicolumn{1}{c|}{69.84}      & 44.94 & \multicolumn{1}{c}{45.06}  \\
\textbf{$D^{PTV}_{30\%}$}               & 36.21 & \multicolumn{1}{c|}{36.23}          & 79.87 & \multicolumn{1}{c|}{80.17}         & 69.40 & \multicolumn{1}{c|}{69.68}      & 44.87 & \multicolumn{1}{c}{44.89}       \\ \hline \hline

\textbf{$D^{Rectum}_{max}$}               & 35.95        & 36.59         & 79.61        & 80.21         & 69.00        & 70.00       & 44.97        & 48.29   \\
\textbf{$D^{Rectum}_{5\%}$}                & 32.87        & 36.23          & 51.42        & 73.36         & 64.97        & 66.09      & 38.95        & 41.79    \\
\textbf{$D^{Rectum}_{mean}$}               & 10.83        & 14.91          & 29.96        & 38.62         & 23.81        & 27.76      & 14.83        & 17.57    \\
\textbf{$D^{Rectum}_{70\%}$}                & 1.68         & 1.69          & 22.87        & 27.33         & 3.11         & 3.06       & 8.59        & 8.61    \\
\end{tabular}
\end{center}
\end{table}

\begin{figure}[htbp]
\begin{center}
\includegraphics[width =0.3 \linewidth]{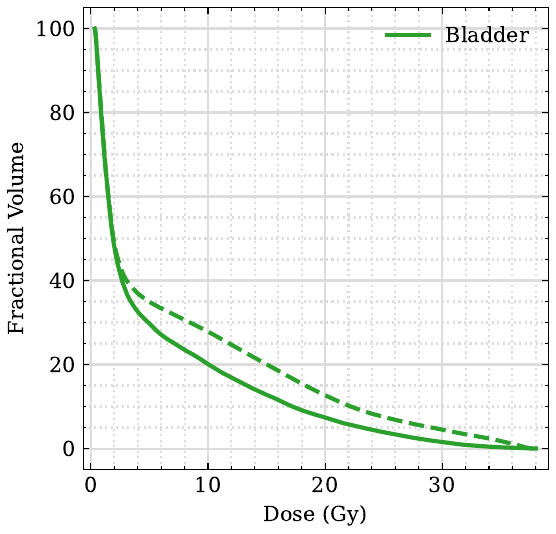}
\includegraphics[width =0.3 \linewidth]{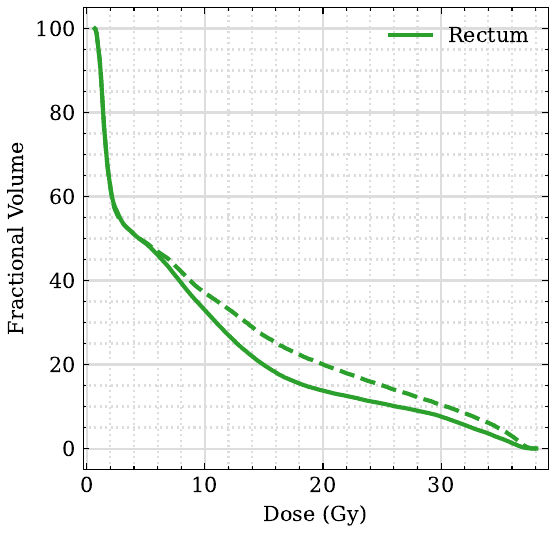}
\includegraphics[width =0.3 \linewidth]{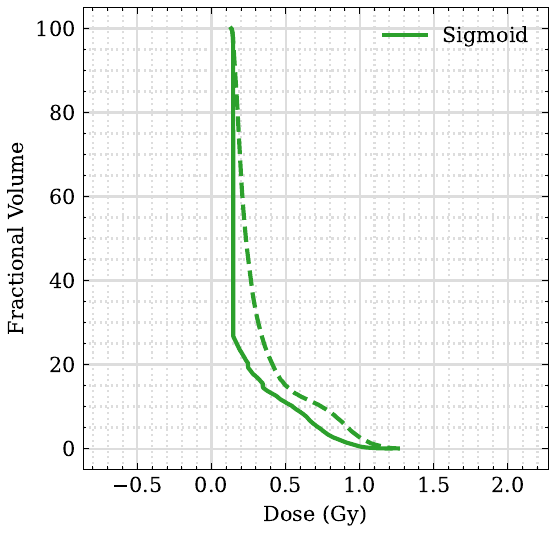}
\includegraphics[width =0.3 \linewidth]{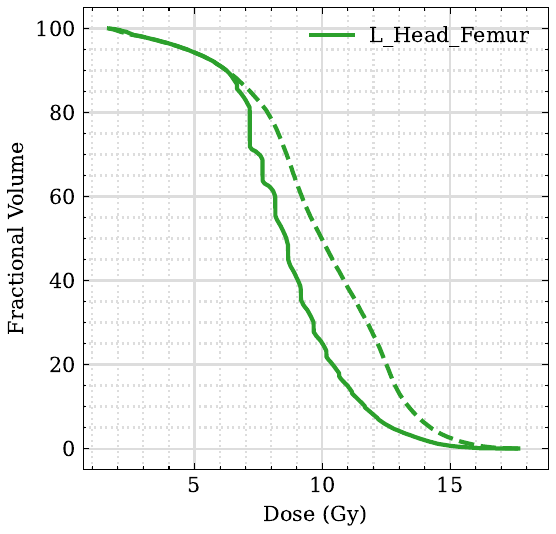}
\includegraphics[width =0.3 \linewidth]{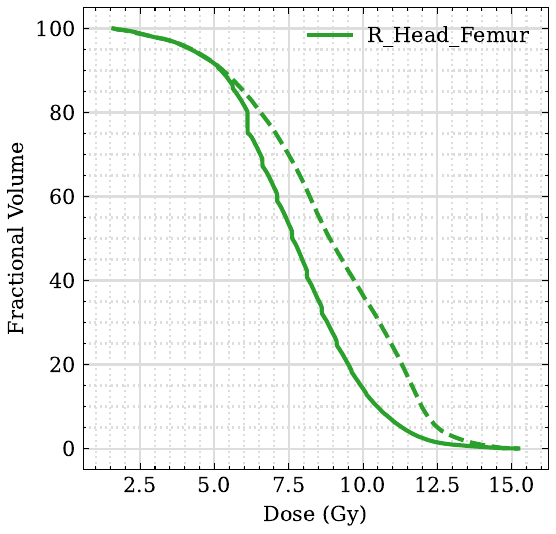}
\caption{\footnotesize Voxel dose differences in histograms for different regions of interest between the plan learned from the inverse learning model and the clinical plan. As shown, plans vary for almost all regions of interest. PTV, although having reduced doses on some voxels, retains similar coverage for 95 percent fractional volume and all other OARs, including sigmoid, remain fairly the same.} \label{Fig:DVH_comparison_all_OARs}
\end{center}
\end{figure}

The results provided so far indicate that improvements are possible for all patients considered in this work in terms of how the mathematical optimization model is capable of guiding the process of finding given solutions to improved solutions. It is worth noting that while improvements are seen for all patients, the level of improving the dose levels for different OARs might be generally subjective as dose levels might already be at optimal or near optimal levels for certain OARs. In the remainder of this section, we provide more detailed results for patient 1 as the representative patient. For patient 1, a detailed description of the objectives used to arrive at an acceptable plan are shown in \cref{Table:Twain_Plan_Objectives}. The table includes the objective types, parameter values for each objective and the weight of the objective in the optimization. Besides the objectives, we extract the dose distribution of the current plan and the masks from the given data. The overlap criteria were also set such that highest priority was given to the PTV. We analyze a prostate cancer patient with a known clinically acceptable plan and recovered objectives from known prediction models. Inputting the known objectives and plan in the $\MIL$ model with the goal of improving DVH objectives on sparing Rectum provides a new plan and a new DVH objective. Our results on a representative prostate case illustrate that $\MIL$ improves the treatment plan and can provide updated objectives. 
Using the detailed overlap criteria, we apply the $\MIL$ model to all major OARs of patient 1 and provide DVH comaprsions. \cref{Fig:DVH_comparison_all_OARs} shows the results of applying ther $\MIL$ model to all major OARs for Patient 1. Our results for a representative patient case show a reduction of 3.3 Gy to 50\% volume of rectum compared to the historical clinical plan, representing a 27.3\% decrease over the average dose to Rectum. The planning target volume (PTV) coverage, while experiencing slightly reduced dose levels, remained comparable and within clinical tolerance in both plans at 95\% volume with only 0.4\% decrease in average dose. The rest of OARs received similar dose levels. Notably, the Rectum DVH curve improves substantially while the PTV coverage remains within clinical targets for the 95\% fractional volume, although insignificantly reduced in dose levels as shown in \cref{Table:Plan_comparisons_patient_1}. Additionally, other structure indicate minimal changes in their DVH curves. Figure 3 compares the plans in voxel level showcasing reduction levels in Rectum and comparable voxel doses in other OARs. Additionally, \cref{Table:Plan_comparisons_patient_1} provides more details on the results of applying $\MIL$ for improving rectum DVH.

\begin{table}[]
\begin{center}
\renewcommand{\arraystretch}{1.5}
\caption{$\MIL$ plan and the original clinical plan for the representative prostate cancer patient in terms of DVH Criteria for select regions of interest including the target (PTV) and key organs at risk (such as rectum and bladder). As shown, Dose-Volume criteria is more desirable for Rectum in the $\MIL$ plan while almost all other regions of interest retain their original behavior. Most notable, 95 percent coverage ($D_{95\%}$) of the target remains intact.}
\label{Table:Plan_comparisons_patient_1}
\begin{tabular}{c|cc|cc|cc|cc}

\textbf{Criteria}               &  \textbf{$D_{max}$}         & \textbf{$D_{max}$}                  & \textbf{$D_{95\%}$}          &  \textbf{$D_{95\%}$}                    &  \textbf{$D_{mean}$} &\textbf{$D_{mean}$}   & \textbf{$D_{30\%}$}         & \textbf{$D_{30\%}$} \\

         & IL               & matRad & IL               & matRad & IL                   & matRad & IL               & matRad   \\
\hline\hline
PTV          & 36.68            & 36.68    & 36.32            & 36.33    & 36.06                & 36.22    & 36.21            & 36.23      \\
Penile Bulb  & 3.29             & 3.04     & 3.21             & 2.91     & 2.05                 & 2.02     & 1.69             & 1.69       \\
Bladder      & 36.65            & 36.65    & 34.22            & 34.35    & 8.28                 & 8.29     & 0.61             & 0.61       \\
\underline{Rectum}       & 35.95            & 36.59    & 32.87            & 36.23    & 10.83                & 14.91    & 1.68             & 1.69       \\
Sigmoid      & 0.82             & 0.82     & 0.35             & 0.35     & 0.03                 & 0.03     & 0.00             & 0.00       \\
L Head Femur & 16.31            & 16.31    & 11.20            & 11.24    & 6.20                 & 6.22     & 2.57             & 2.55       \\
R Head Femur & 15.61            & 15.61    & 10.54            & 10.55    & 6.50                 & 6.52     & 3.42             & 3.25 \\
\hline
\end{tabular}
\end{center}
\end{table}

The analysis over the representative patient shows that the MIL approach is capable of proposing new realizations of $\IPP$ along with improved plans for such realizations. The results indicate that using inverse learning methods can improve the DV objectives in radiation therapy treatment plans. The outcomes illustrate the potential for data-driven frameworks capable of identifying the best improvements for groups of patients with similar attributes.

\subsection{Pareto Frontier Implications} \label{sec:results_Pareto}

In addition to single OAR comparisons of existing plans and $\MIL$ improved plans, we conduct a series of Pareto Frontier comparisons across the two sets of objectives to showcase how using $\MIL$ results in different solutions that might be inaccessible to the original treatment plan objective set. For this analysis, we employ the open-access software matRad and optimize the radiation therapy treatment planning problem using the initial set of parameters and the improved set of parameters for the maximum DVH objective of the 20\% fractional volume of the rectum. Then, four Pareto optimal plans are generated by varying the weights of the objectives for the maximum DVH objective of the 20\% fractional volume of the rectum and the minimum DVH of the 95\% fractional volume of the PTV. 

\begin{figure}[htbp]
\begin{center}
\includegraphics[width =1 \linewidth]{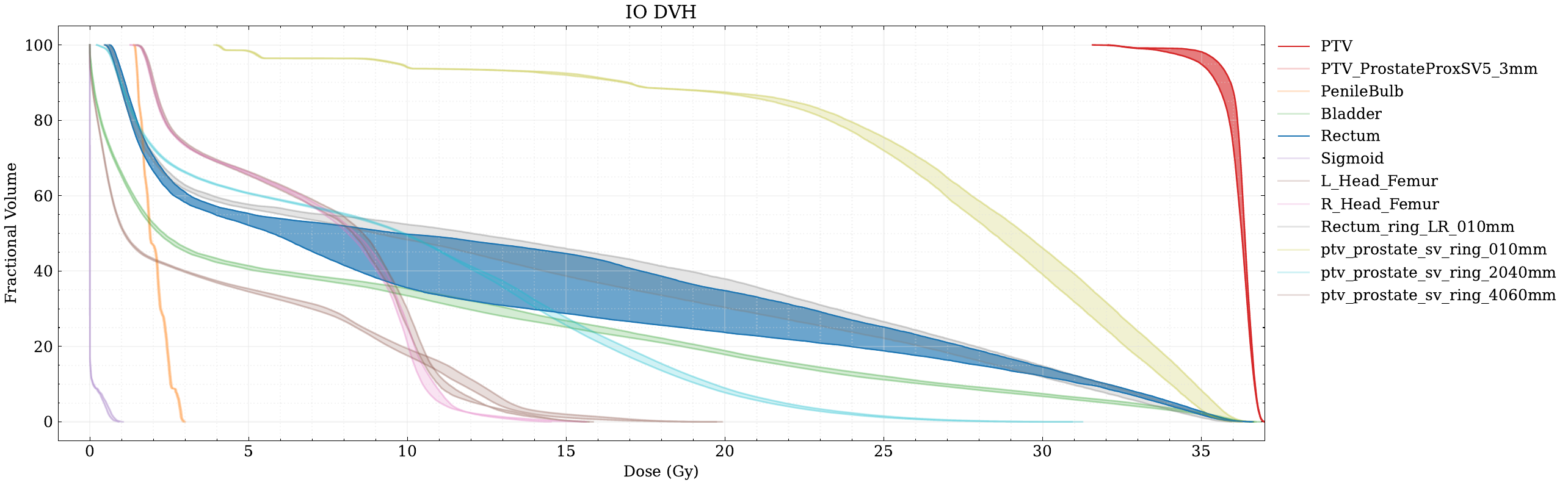}
\includegraphics[width =1 \linewidth]{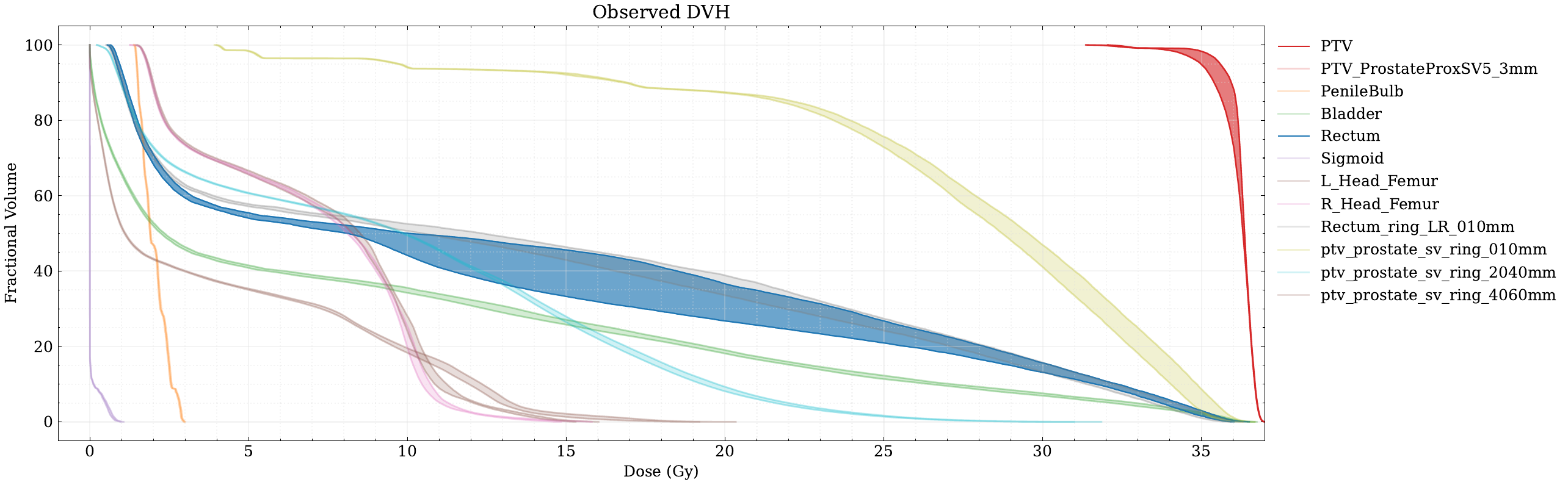}
\caption{\footnotesize Dose-volume histogram (DVH) comparisons for Pareto optimal plans generated using the original (bottom) and MIL-optimized (top) objective parameters. Each Pareto optimal plan represents a different trade-off between the competing objectives of rectum sparing (20\% fractional volume) and PTV coverage (95\% fractional volume). The MIL-optimized objective parameters result in a Pareto frontier that consistently achieves better rectum sparing across Pareto optimal plans, as evidenced by the lower rectum DVH curves, while maintaining similar PTV coverage. This demonstrates the effectiveness of the MIL approach in exploring the trade-off space and uncovering previously unattainable solutions that improve organ-at-risk sparing without compromising target dose delivery.} \label{Fig:Pareto Comparisons}
\end{center}
\end{figure}

\begin{figure}[htbp]
\begin{center}
\includegraphics[width =0.18 \linewidth]{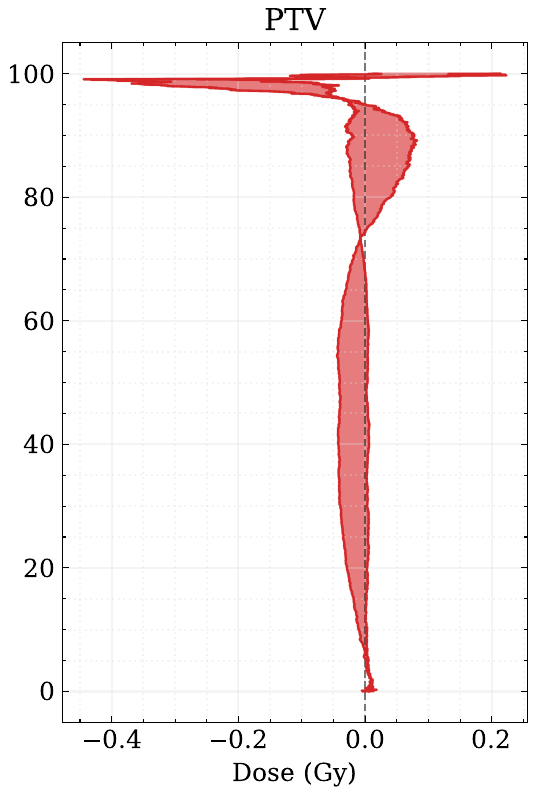}
\includegraphics[width =0.18 \linewidth]{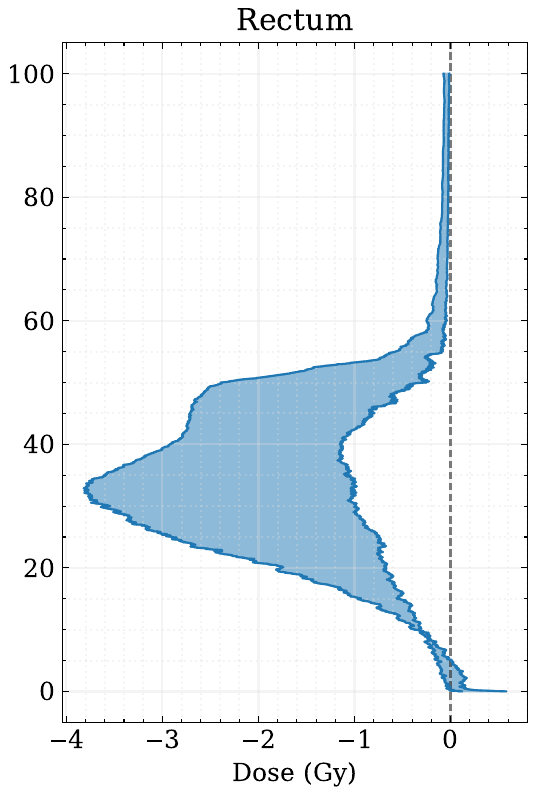}
\includegraphics[width =0.18 \linewidth]{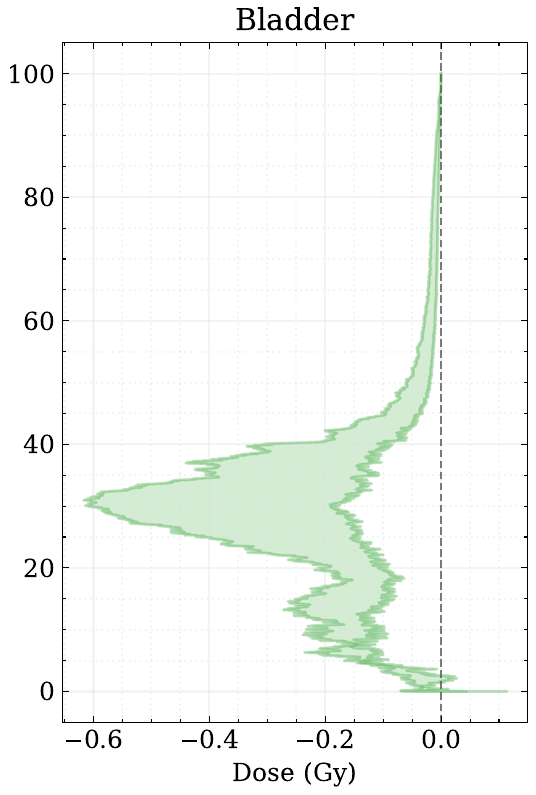}
\includegraphics[width =0.18 \linewidth]{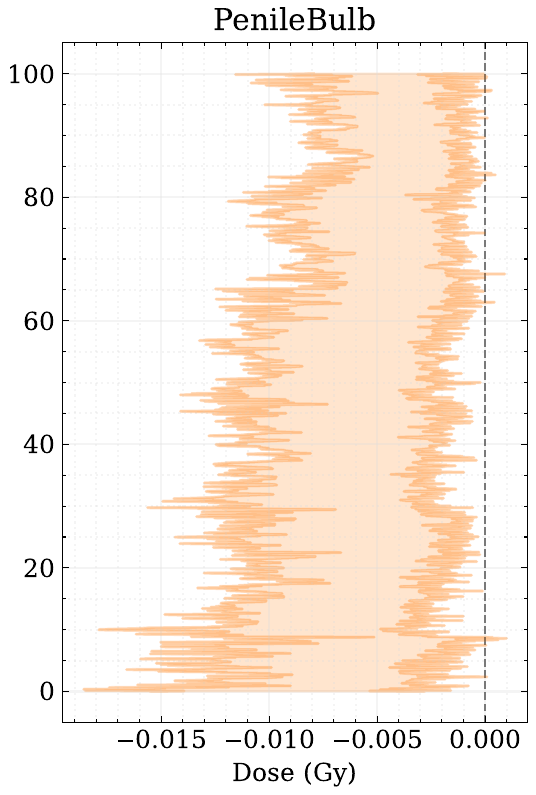}
\includegraphics[width =0.18 \linewidth]{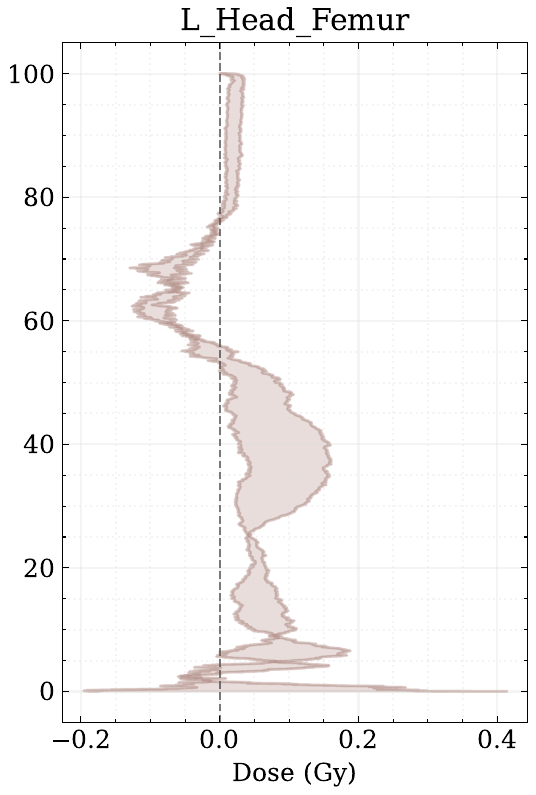}
\includegraphics[width =0.18 \linewidth]{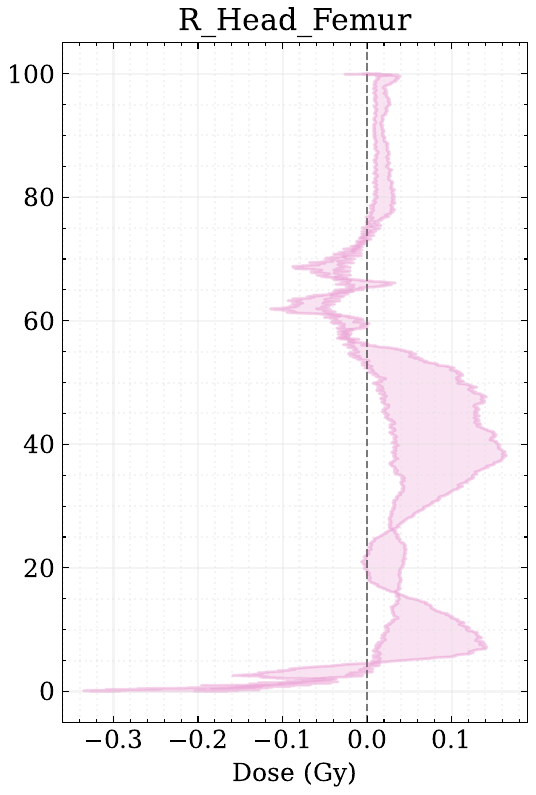}
\includegraphics[width =0.18 \linewidth]{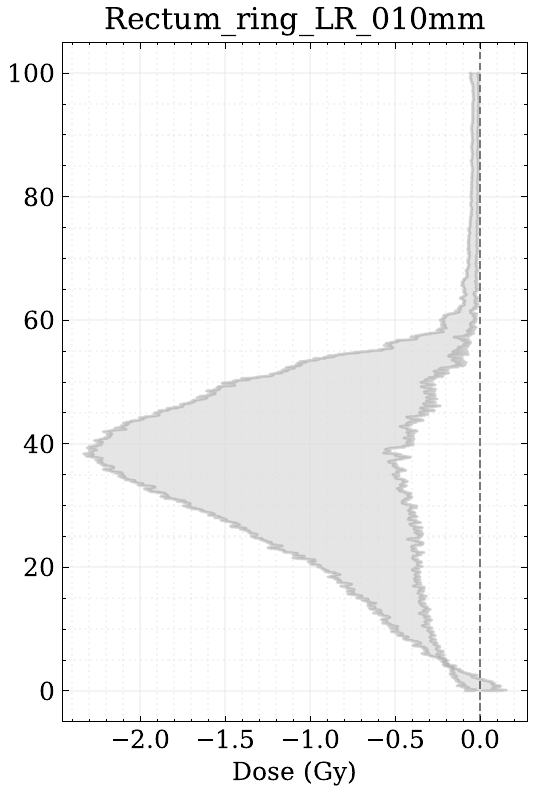}
\includegraphics[width =0.18 \linewidth]{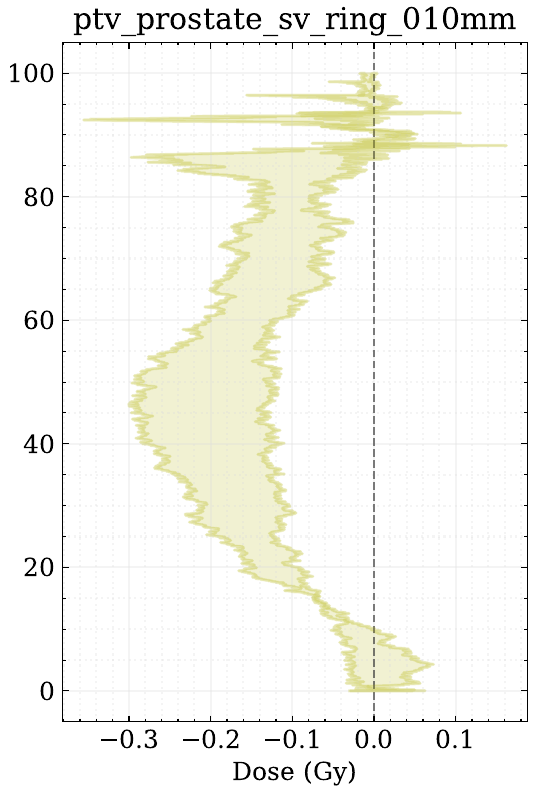}
\includegraphics[width =0.18 \linewidth]{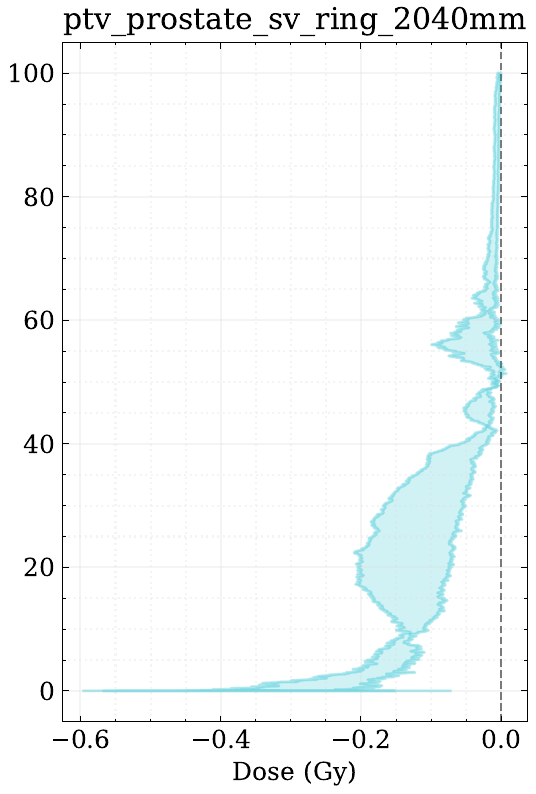}
\includegraphics[width =0.18 \linewidth]{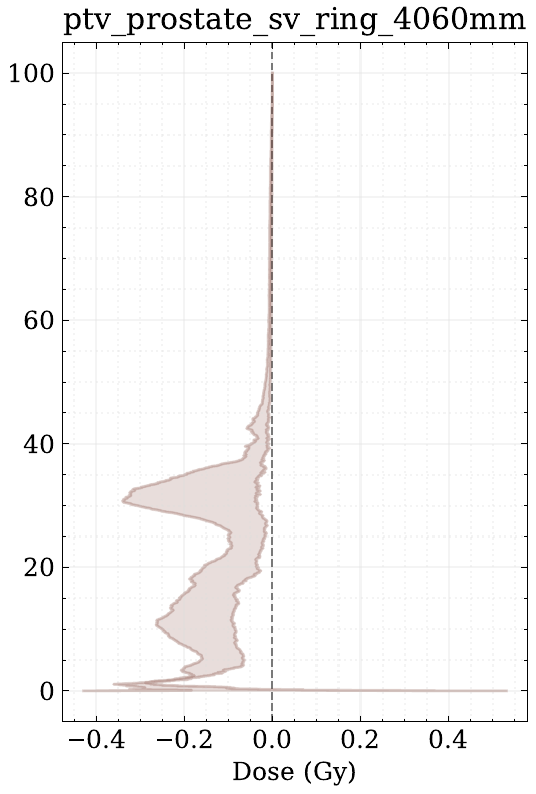}
\caption{\footnotesize \footnotesize Detailed comparisons for each OAR and target structure, across the Pareto optimal plans generated using the original  and MIL-optimized objective parameters in comparison to the observed plan (vertical dashed line). The rectum DVHs exhibit the most notable improvements, with consistent dose reductions across all Pareto optimal plans. The PTV and other OARs, such as the bladder, penile bulb, and femoral heads, maintain similar DVH curves between the two sets of plans, indicating that the MIL-driven rectum sparing improvements are achieved without compromising dose distributions in other critical structures. These OAR-specific comparisons highlight the MIL model's ability to fine-tune trade-offs and optimize dose-volume objectives for individual structures while preserving overall plan quality.} \label{Fig:Pareto Comparisons_OARS}
\end{center}
\end{figure}

As shown in \cref{Fig:Pareto Comparisons}, by updating the parameters of the dose limit for the maximum DVH objective of the 20\% fractional volume of the rectum, a superior DVH figure is achieved that does not compromise the dose level on the PTV but reduces shifts the overall rectum curve towards less dose. \cref{Fig:Pareto Comparisons_OARS} showcases these improvements over Pareto solutions in greater detail.

\section{Discussions and Conclusions} \label{sec:discussions}

In this study, we developed a novel inverse optimization-based framework for improving radiation therapy treatment plans using observed clinical plans. The proposed methodology, termed Modified Inverse Learning (MIL), allows for the automated discovery of improved treatment plans by optimizing the parameters of the treatment planning optimization problem while preserving the preferences of the original optimizer. The MIL approach was applied to the specific case of prostate cancer radiation therapy, where the rectum is often the primary dose-limiting organ at risk (OAR) \citep{wang2013quality}.

Our results demonstrate the potential of the MIL framework to generate improved treatment plans that reduce the dose to critical structures, such as the rectum, while maintaining target coverage. In the representative prostate cancer case, the MIL plan achieved a substantial reduction in rectum dose, with a 27.3\% decrease in the average dose and a 3.3 Gy reduction in the dose to 50\% of the rectum volume compared to the clinical plan. Notably, these improvements were achieved without compromising the dose coverage of the planning target volume (PTV), which remained within clinical tolerances. The MIL framework's ability to optimize treatment plans while preserving the original optimizer's preferences is a significant advantage over existing inverse optimization approaches \citep{chan2020inverseconstraint}.

The MIL framework's performance was further validated on a cohort of four prostate cancer patients, demonstrating consistent improvements in rectum sparing across all cases. The results highlight the MIL approach's robustness and its potential to generate high-quality treatment plans that are tailored to individual patient anatomy and clinical goals. By automating the plan optimization process and incorporating prior knowledge from clinical plans, the MIL framework offers a valuable tool for enhancing the efficiency and quality of radiation therapy treatment planning.

The exploration of the Pareto frontier for the multi-objective optimization problem revealed valuable insights into the trade-offs between target coverage and OAR sparing. By comparing the Pareto frontiers generated using the original and MIL-optimized objective parameters, we demonstrated that the MIL approach could lead to the discovery of previously unattainable solutions that improve OAR sparing without compromising target coverage. This finding underscores the importance of parameter optimization in radiation therapy treatment planning and highlights the MIL framework's potential to navigate complex trade-off spaces effectively.

The MIL framework's ability to handle multiple OARs and various types of dose-volume constraints demonstrates its flexibility and applicability to a wide range of radiation therapy treatment planning scenarios. Although the present study focused on prostate cancer, the MIL approach can be easily adapted to other cancer sites and treatment modalities, such as head and neck cancer or stereotactic body radiation therapy (SBRT). Further research is warranted to investigate the performance of the MIL framework in these diverse clinical settings and to assess its potential impact on patient outcomes and treatment efficiency.

One limitation of the current study is the assumption that the feasible set of the treatment planning optimization problem is convex and can be accurately approximated using linear constraints. While this assumption simplifies the optimization problem and facilitates the interpretation of the results, it may not fully capture the complexities of the actual feasible set. Future work should explore the incorporation of more sophisticated dose calculation algorithms and non-linear constraints to improve the realism and accuracy of the MIL framework.

Another avenue for future research is the integration of the MIL framework with data-driven approaches for predicting achievable dose-volume objectives \citep{mahmoudzadeh2015robust}. By combining the strengths of inverse optimization and machine learning, it may be possible to develop a comprehensive treatment planning system that leverages prior clinical knowledge, patient-specific anatomical information, and optimization algorithms to generate high-quality, personalized treatment plans.

In conclusion, the Modified Inverse Learning (MIL) framework presented in this study offers a powerful and flexible approach for optimizing radiation therapy treatment plans using observed clinical plans. By automating the discovery of improved plans while preserving the preferences of the original optimizer, the MIL approach has the potential to significantly enhance the quality and efficiency of radiation therapy treatment planning. The promising results obtained for prostate cancer patients motivate further research into the application of the MIL framework to other cancer sites and treatment modalities, as well as its integration with data-driven approaches for personalized treatment planning.

\appendix
\section{Proofs of Statements}
\label{app1}

\newproof{pop1}{Proof of Proposition \ref{ILG_feasible}}
\begin{pop1} 
Let $\bx \in \Omega$ be a point on the boundary of $\Omega$ such that $\hat{\ba}_i \bx = \hat{b}^0_i$ for some $i \in \hat{\cJ}$. We show that $(\hat{\ba}_i, \bx, \left [e_i; \,0 \right], \hat{\bb}_0)$ is feasible for $\IL_g$. First, note that $\bx \in \Omega$ implies that $\bx$ satisfies constraints \eqref{MILPrimal Feasblility1} and \eqref{MILPrimal Feasblility2} in $\IL_g$.  Next, we show that $(\hat{\ba}_i, \bx, \left [e_i; \,0 \right], \hat{\bb}_0)$ satisfies constraints \eqref{MILDualFeas1} - \eqref{MILDualFeas2}. Since $e_i$ is the $i^{th}$ basis vector, we have:
\begin{align*}
    \bA' \left [e_i; \,0 \right] &= \hat{\ba}_i
\end{align*}
which satisfies constraint \eqref{MILDualFeas1}. Furthermore,
\begin{align*}
    \hat{\ba}_i' \bx &= \hat{b}^0_i \\
    &= \left [\hat{\bb}_0; \,\bar{\bb} \right]' \left [e_i; \,0 \right]
\end{align*}
satisfying constraint \eqref{MILStrongDual}. Finally, since $e_i \geq 0$ and $\sum_{j \in \hat{\cJ}} e_{ij} =1$, constraints \eqref{MILRegularization} and \eqref{MILDualFeas2} are also satisfied. Therefore, $(\hat{\ba}_i, \bx, \left [e_i; \,0 \right], \hat{\bb}_0)$ satisfies all constraints of $\IL_g$ and is feasible for $\IL_g$. 
\end{pop1}

\newproof{pop2}{Proof of Proposition \ref{prop:ILG_feasiblity_properties}}
\begin{pop2}
We prove each part separately:

\begin{enumerate}
    \item Let $(\bc, \bx, \by, \hat{\bb}) \in \Phi (\hat{\cJ})$. Then, $(\bc, \bx, \by, \hat{\bb})$ satisfies constraints \eqref{MILPrimal Feasblility1} - \eqref{MILrange} of $\IL_g$. In particular, $\bx$ satisfies:
    \begin{align*}
        \bar{\bA} \bx &\leq \bar{\bb} \\
        \hat{\bA} \bx &\leq \hat{\bb} \\
        \hat{\bb}_L &\leq \hat{\bb} \leq \hat{\bb}_U
    \end{align*}
    Combining the last two inequalities, we have:
    \begin{align*}
        \hat{\bb}_L &\leq \hat{\bA} \bx \leq \hat{\bb}_U
    \end{align*}
    then $\bx \in \left\{ \bx \in \mathbb{R}^n | \bar{\bA} \bx \leq \bar{\bb}, \hat{\bb}_L\leq\hat{\bA} \bx\leq \hat{\bb}_U \right\}$ and $\Phi (\hat{\cJ})_{\bx} \subseteq \left\{ \bx \in \mathbb{R}^n | \bar{\bA} \bx \leq \bar{\bb}, \hat{\bb}_L\leq\hat{\bA} \bx\leq \hat{\bb}_U \right\}$.
    
    Conversely, let $\bx \in \left\{ \bx \in \mathbb{R}^n | \bar{\bA} \bx \leq \bar{\bb}, \hat{\bb}_L\leq\hat{\bA} \bx\leq \hat{\bb}_U \right\}$. Then, there exists $\hat{\bb}$ such that $\hat{\bb}_L \leq \hat{\bb} \leq \hat{\bb}_U$ and $\hat{\bA} \bx \leq \hat{\bb}$. By Proposition \ref{ILG_feasible}, there exist $\bc$ and $\by$ such that $(\bc, \bx, \by, \hat{\bb})$ is feasible for $\IL_g$. Therefore, $\bx \in \Phi (\hat{\cJ})_{\bx}$, and hence $\left\{ \bx \in \mathbb{R}^n | \bar{\bA} \bx \leq \bar{\bb}, \hat{\bb}_L\leq\hat{\bA} \bx\leq \hat{\bb}_U \right\} \subseteq \Phi (\hat{\cJ})_{\bx}$.
    
    \item Let $\hat{\cJ}_1$ and $\hat{\cJ}_2$ be such that $\hat{\cJ}_2 \subseteq \hat{\cJ}_1$. Let $(\bc, \bx, \by, \hat{\bb}) \in \Phi (\hat{\cJ}_2)$. Then, $(\bc, \bx, \by, \hat{\bb})$ satisfies constraints \eqref{MILPrimal Feasblility1} - \eqref{MILrange} of $\IL_g$ with $\hat{\cJ} = \hat{\cJ}_2$. Since $\hat{\cJ}_2 \subseteq \hat{\cJ}_1$, these constraints are a subset of the constraints of $\IL_g$ with $\hat{\cJ} = \hat{\cJ}_1$. Therefore, $(\bc, \bx, \by, \hat{\bb})$ also satisfies constraints \eqref{MILPrimal Feasblility1} - \eqref{MILrange} of $\IL_g$ with $\hat{\cJ} = \hat{\cJ}_1$, and hence $(\bc, \bx, \by, \hat{\bb}) \in \Phi (\hat{\cJ}_1)$. This shows that $\Phi (\hat{\cJ}_2) \subseteq \Phi (\hat{\cJ}_1)$.
    
    Finally, note that $\Phi$ is the feasible set of $\IL_g$ with $\hat{\cJ} = \emptyset$. Since $\emptyset \subseteq \hat{\cJ}$ for all $\hat{\cJ}$, we have $\Phi \subseteq \Phi (\hat{\cJ})$ for all $\hat{\cJ}$.
\end{enumerate}
\end{pop2}

\newproof{pop3}{Proof of Theorem \ref{thm:ILG_solution}}
\begin{pop3}
Let $(\bc, \bx, \by, \hat{\bb})$ be an optimal solution for $\IL_g({\bx_0, \bA, \bar{\bb}, \omega})$. By the strong duality constraint \eqref{MILStrongDual}, we have:
\begin{align}
    \bc' \bx &= \left [\hat{\bb}; \,\bar{\bb} \right]' \left [\hat{\by}; \,\bar{\by} \right] \label{eq:thm1_1}
\end{align}

Since $\bx_0 \in \Phi(\hat{\cJ})$, there exist $\hat{\bb}_0$ and $\by_0$ such that $(\bc, \bx_0, \by_0, \hat{\bb}_0)$ is feasible for $\IL_g({\bx_0, \bA, \bar{\bb}, \omega})$. Therefore, we have:
\begin{align}
    \bc' \bx_0 &= \left [\hat{\bb}_0; \,\bar{\bb} \right]' \left [\hat{\by}_0; \,\bar{\by}_0 \right] \label{eq:thm1_2}
\end{align}

Now, suppose for contradiction that $\bc' \bx < \bc' \bx_0$. Then, from \eqref{eq:thm1_1} and \eqref{eq:thm1_2}, we have:
\begin{align*}
    \left [\hat{\bb}; \,\bar{\bb} \right]' \left [\hat{\by}; \,\bar{\by} \right] &< \left [\hat{\bb}_0; \,\bar{\bb} \right]' \left [\hat{\by}_0; \,\bar{\by}_0 \right]
\end{align*}

Since $\hat{\by}, \hat{\by}_0 \geq 0$ and $\sum_{j \in \hat{\cJ}} \hat{y}_j = \sum_{j \in \hat{\cJ}} \hat{y}_{0j} = 1$, we have:
\begin{align}
    \hat{\bb}' \hat{\by} &< \hat{\bb}_0' \hat{\by}_0 \label{eq:thm1_3}
\end{align}

Consider the solution $(\bc, \bx_0, \by_0, \hat{\bb})$. This solution is feasible for $\IL_g({\bx_0, \bA, \bar{\bb}, \omega})$ because $\bx_0 \in \Phi(\hat{\cJ})$ and $\hat{\bb}_L \leq \hat{\bb} \leq \hat{\bb}_U$. Moreover, from \eqref{eq:thm1_3}, we have:
\begin{align*}
    \omega\cD(\bx_0, \bx_0) - (1 - \omega) (\left [\hat{\bb}; \,\bar{\bb} \right]' \left [\hat{\by}_0; \,\bar{\by}_0 \right]) &< \omega\cD(\bx, \bx_0) - (1 - \omega) (\left [\hat{\bb}; \,\bar{\bb} \right]' \left [\hat{\by}; \,\bar{\by} \right])
\end{align*}
which contradicts the optimality of $(\bc, \bx, \by, \hat{\bb})$. Therefore, we must have $\bc' \bx \geq \bc' \bx_0$.
\end{pop3}

\newproof{pop4}{Proof of Proposition \ref{prop:ILg_Ab}}
\begin{pop4}
We prove each part separately:

\begin{enumerate}
    \item Let $\omega = 1$. Consider the solution $(\hat{\ba}_i, \bx_0, \left [e_i; \,0 \right], \hat{\bb}_0)$ where $i \in \hat{\cJ}$ is such that $\hat{\ba}_i \bx_0 = \hat{b}^0_i$. By Proposition \ref{ILG_feasible}, this solution is feasible for $\IL_g({\bx_0, \bA, \bar{\bb}, 1})$. Moreover, for any feasible solution $(\bc, \bx, \by, \hat{\bb})$, we have:
    \begin{align*}
        \cD(\bx, \bx_0) - 0 &\geq \cD(\bx_0, \bx_0) - 0 = 0
    \end{align*}
    Therefore, $(\hat{\ba}_i, \bx_0, \left [e_i; \,0 \right], \hat{\bb}_0)$ is an optimal solution for $\IL_g({\bx_0, \bA, \bar{\bb}, 1})$ with $\bx^* = \bx_0$.
    
    \item Let $\omega = 0$. Consider the solution $(\hat{\ba}_i, \hat{\bx}, \left [e_i; \,0 \right], \hat{b} e_i)$ where $i \in \hat{\cJ}$, $\hat{\bx} \in \mathbb{R}^n$ is such that $\hat{\ba}_i \hat{\bx} = \hat{b}$ for some $\hat{b} \in \{\hat{\bb}_L\} \cup \{\hat{\bb}_U\}$. By Proposition \ref{ILG_feasible}, this solution is feasible for $\IL_g({\bx_0, \bA, \bar{\bb}, 0})$. Moreover, for any feasible solution $(\bc, \bx, \by, \hat{\bb})$, we have:
    \begin{align*}
        0 - \left [\hat{\bb}; \,0 \right]' \left [\hat{\by}; \,0 \right] &= - \hat{b} \leq - \hat{b}_k
    \end{align*}
    for some $k \in \{1,\ldots,|\hat{\cJ}|\}$ and $\hat{b}_k \in \{\hat{\bb}_L\} \cup \{\hat{\bb}_U\}$. Therefore, $(\hat{\ba}_i, \hat{\bx}, \left [e_i; \,0 \right], \hat{b} e_i)$ is an optimal solution for $\IL_g({\bx_0, \bA, \bar{\bb}, 0})$ with $\hat{\ba}_i \bx^* = \hat{b}$ for some $\hat{b} \in \{\hat{\bb}_L\} \cup \{\hat{\bb}_U\}$.
\end{enumerate}
\end{pop4}

\newproof{pop5}{Proof of Theorem \ref{thm:MIL_ILG}}
\begin{pop5}
Let $\{\bx, \hat{b}^* \} $ be an optimal solution for $\MIL$ with $\hat{\cJ} = \{j\}$. We construct a solution $\{\bc^*, \bx^*, \by^*, \hat{\bb}^* \}$ for $\IL_g({\bx_0, \bA, \bar{\bb}, \omega})$ as follows:

1. Set $\bx^* = \bx$ from the optimal solution of $\MIL$.
2. Set $\hat{\bb}^* = \hat{b}^* e_j$ where $e_j$ is the $j$-th standard basis vector.
3. Set $\by^* = \left [e_j; \,0 \right]$.
4. Set $\bc^* = \bA' \by^*$.

We now show that $\{\bc^*, \bx^*, \by^*, \hat{\bb}^* \}$ is feasible and optimal for $\IL_g({\bx_0, \bA, \bar{\bb}, \omega})$.

Feasibility: Since $\{\bx^*, \hat{b}^* \}$ is feasible for $\MIL$, we have:
\begin{align*}
    \bar{\bA} \bx^* &\leq \bar{\bb} \\
    \hat{\ba}_j \bx^* &= \hat{b}^*
\end{align*}
Therefore, $\bx^*$ satisfies constraints \eqref{MILPrimal Feasblility1} and \eqref{MILPrimal Feasblility2} of $\IL_g$. Moreover, by construction, $\by^*$ and $\bc^*$ satisfy constraints \eqref{MILDualFeas1}, \eqref{MILRegularization}, and \eqref{MILDualFeas2}. Finally, we have:
\begin{align*}
    (\bc^*)' \bx^* &= (\bA' \by^*)' \bx^* \\
                  &= (\by^*)' \bA \bx^* \\
                  &= \left [e_j; \,0 \right]' \left [\hat{b}^*; \,\bar{\bb} \right] \\
                  &= \hat{b}^* \\
                  &= \left [\hat{\bb}^*; \,\bar{\bb} \right]' \left [\hat{\by}^*; \,\bar{\by}^* \right]
\end{align*}
which satisfies constraint \eqref{MILStrongDual}. Therefore, $\{\bc^*, \bx^*, \by^*, \hat{\bb}^* \}$ is feasible for $\IL_g({\bx_0, \bA, \bar{\bb}, \omega})$.

Optimality: Suppose for contradiction that $\{\bc^*, \bx^*, \by^*, \hat{\bb}^* \}$ is not optimal for $\IL_g({\bx_0, \bA, \bar{\bb}, \omega})$. Then, there exists a feasible solution $\{\tilde{\bc}, \tilde{\bx}, \tilde{\by}, \tilde{\bb} \}$ such that:
\begin{align*}
    \omega\cD(\tilde{\bx}, \bx_0) - (1 - \omega) (\left [\tilde{\bb}; \,\bar{\bb} \right]' \left [\tilde{\by}; \,\bar{\by} \right]) &< \omega\cD(\bx^*, \bx_0) - (1 - \omega) (\left [\hat{\bb}^*; \,\bar{\bb} \right]' \left [\hat{\by}^*; \,\bar{\by}^* \right]) \\
    &= \omega\cD(\bx^*, \bx_0) - (1 - \omega) \hat{b}^*
\end{align*}
Since $\tilde{\by} \geq 0$ and $\sum_{j \in \hat{\cJ}} \tilde{y}_j = 1$, we have:
\begin{align*}
    \omega\cD(\tilde{\bx}, \bx_0) - (1 - \omega) \tilde{b}_j &< \omega\cD(\bx^*, \bx_0) - (1 - \omega) \hat{b}^*
\end{align*}
where $\tilde{b}_j = \hat{\ba}_j \tilde{\bx}$. But this contradicts the optimality of $\{\bx^*, \hat{b}^* \}$ for $\MIL$. Therefore, $\{\bc^*, \bx^*, \by^*, \hat{\bb}^* \}$ must be optimal for $\IL_g({\bx_0, \bA, \bar{\bb}, \omega})$.
\end{pop5}

\bibliographystyle{elsarticle-harv} 
\bibliography{MLIORT}



\end{document}